\documentclass[10pt]{article}
\usepackage{amsmath, amsthm, amsfonts}
\usepackage{subcaption}
\usepackage{graphicx}

\usepackage{vmargin}

\theoremstyle{definition}

\theoremstyle{remark}

\newtheorem{teorema}{Theorem}[section]
\newtheorem{propo}[teorema]{Proposition}

\newtheorem{lema}[teorema]{Lemma}
\newtheorem{obs}[teorema]{Observation}

\newtheorem{corolario}[teorema]{Corollary}

\title{An inverse problem for a linear system of dispersive equations}
\author{Deissy Marcela Pizo\footnote{Universidad del Valle, Departamento de Matem\'aticas, Cali, Colombia, email: deissy.pizo@correounivalle.edu.co }, Juan Carlos Mu\~noz Grajales \footnote{Universidad del Valle, Departamento de Matem\'aticas, Cali, Colombia, Corresponding author, email: juan.munoz@correounivalle.edu.co } } 

\begin{document}
\maketitle

\begin{abstract}
This paper addresses the inverse problem of identifying the linear velocity coefficient in a linear system governed by two Benjamin-Bona-Mahony-type equations, which model the displacement of water waves propagating along the surface of a shallow channel, incorporating effects of dispersion and topography. To solve this, we reformulate the inverse problem as a restricted minimization problem (RMP), aimed at optimizing a suitably regularized objective functional. We use numerical techniques, specifically the iterative L-BFGS-B algorithm implemented in the Dolfin-Adjoint-Python-SciPy libraries, to solve the RMP effectively. Following methodologies similar to those in \cite{JCPipicanoSosa}, we establish a local stability result for the RMP. Additionally, through numerical simulations, we demonstrate the effectiveness of the proposed identification method in determining the linear velocity coefficient in Boussinesq-type systems.

\end{abstract}

\vspace{0.5cm}


\vspace{0.5cm}


\noindent Secondary: 65N40, 35Q53

\vspace{0.5cm}

\section{Introduction}

Inverse problems arise in diverse fields, including geophysics, engineering, medical imaging, and computer vision. In geophysics, they enable the analysis of seismic data to infer subsurface properties, such as mineral distributions and 
locations of oil and gas deposits. In engineering, inverse problems help identify structural defects, optimize system designs, and monitor machinery health. In medical imaging, they are used to reconstruct internal body structures from data obtained by devices like X-rays, MRI, and CT scans. Likewise, in computer vision, inverse problems allow for the recovery of object shapes and scene layouts from images or videos.

Inverse problems arise when deducing unknown properties of physical systems from measurements or observations. They are particularly relevant in the study of Partial Differential Equations (PDEs), where unknown parameters or functions must be inferred from solution measurements \cite{Ding}, \cite{Liu1}, \cite{Liu2}. For example, medical imaging often involves estimating contrast agent distributions from X-ray data, while in geophysics, subsurface rock compositions are inferred from seismic readings. Addressing these challenges requires developing mathematical models that relate unknowns to observed data and using these models to estimate hidden quantities. This process demands a deep understanding of the underlying PDEs and their properties, along with advanced tools from optimization theory, numerical analysis, and statistical inference. Consequently, inverse problems represent a dynamic research field with far-reaching implications across science and engineering, driving advancements in numerous disciplines.

In \cite{Grajales2004} is analyzed well-posedness of the Boussinesq-type system
\begin{equation}\label{KdVe}
  \begin{array}{c l}
   M(\xi) \eta_t+\partial _\xi \left[\left(1+\frac{ \tilde{\alpha} \eta}{M\left(\xi\right)}\right)u\right]-\frac{\beta}{6}\partial^2 _\xi \left(M(\xi) \eta_t\right)&=0, \:\: \left(\xi,t\right)\in \mathbb{R}\times \left[0,\infty\right)\\
  u_t+\eta_\xi+\frac{ \tilde{\alpha} }{2}\partial _\xi\left[\left(\frac{u}{M\left(\xi\right)}\right)^2\right]-\frac{\beta}{6}\partial^2 _\xi \left(u_t\right)&=0, \:\: \left(\xi,t\right)\in \mathbb{R}\times \left[0,\infty\right).
   \end{array}
\end{equation}
These equations describe the bidirectional propagation of long waves with small amplitude over a shallow uni-dimensional channel with a variable depth, which may be highly oscillatory, i.e. when the variation scale of the depth irregularities is small with respect to the characteristic wavelength. In \cite{Grajales2004} was performed a study on global existence, regularity and continuous dependence on initial conditions of an initial-value problem associated to system \eqref{KdVe}, where it is considered $0< \beta<<1$ (weakly dispersive regime) and $0 <\tilde{\alpha}<<1$ (weakly nonlinear regime). We point out that, within this parameter regime, the Boussinesq system \eqref{KdVe} is a valid asymptotic approximation of the full nonlinear potential theory equations. This is the reason for the special interest on this physical regime.  The function $u\left(\xi,t\right)$ represents the fluid velocity at a fixed depth, $\eta\left(\xi,t\right)$  denotes the wave elevation and $M\left(\xi\right)$ is a coefficient which depends on the submerged topography of the channel. In the case of a constant depth, the coefficient $M\left(\xi\right)$ is identically one,
and the Boussinesq system \eqref{KdVe} reduces to the one considered by Bona and Chen \cite{Bona1998}. 

As an initial step in investigating inverse problems associated to dispersive-type systems, in the present work we consider from the analytical and numerical point of view the limit linear regime of system in \eqref{KdVe}, i.e. taking $\tilde{\alpha} =0$, proposed on a finite interval with length $L$.
Using the change of variable $N=M\eta$, $V=u$ and introducing the variable coefficient $c(\xi)=\frac{1}{M(\xi)}$ into system \eqref{KdVe} with $\tilde{\alpha}=0$, we obtain the Boussinesq-type system given by 
\begin{equation}\label{KdVe1}
  \begin{array}{c l}
   N_t+V_{\xi}-\frac{\beta}{6}N_{\xi\xi t}&=0 \:\:\: ,\:\:\left(\xi,t\right)\in \left[0,L\right]\times\left[0,T\right].\\
   V_t+(c(\xi) N)_{\xi}-\frac{\beta}{6}V_{\xi\xi t}&=0,
   \end{array}
\end{equation}
subject to initial conditions 
\begin{equation}\label{datosiniciales1}
\begin{array}{c l}
 N(\xi,0)&=N_0(\xi),\\
 V(\xi,0)&=V_0(\xi),\:\:\: \xi \in \left[0,L\right],
\end{array}
\end{equation}
and Dirichlet boundary-value conditions
\begin{equation}\label{condfrontera1}
\begin{array}{c l}
 N(0,t)=N(L,t)=0,\:\: t\in \left[0,T\right]\\
 V(0,t)=V(L,t)=0,\:\: t\in \left[0,T\right],
\end{array}
\end{equation}
where $N(\xi,t), V(\xi,t) $ are related to the wavelength and fluid velocity, respectively, and the coefficient $c(\xi)$ depends only on the spatial variable $\xi$ and it is related to the channel topography. We point out that this change of variables is introduced only to make easier the establishment of the theoretical results. 

This paper has two primary objectives. Firstly, we analytically explore the inverse problem of reconstructing the linear long wave speed, denoted as $c(\xi)$, within the linear system \eqref{KdVe1}. Our approach involves leveraging single measurements $N_T(\xi)$, and $V_T(\xi)$ of the quantities $N(\xi,t), V(\xi,t)$, respectively, at the final time $t=T$. We propose a method rooted in least-squares minimization, augmented by a Tikhonov regularization term added to the objective functional. The reader is encouraged to look at the original papers by Tikhonov \cite{Tikhonov, Tikhonov19632, Tikhonov19631} for the basic results on this technique.

Secondly, we aim to devise a numerical technique for approximating the solution of the aforementioned inverse problem within the full nonlinear system \eqref{KdVe}. This involves solving a restricted minimization problem (RMP) associated with the Tikhonov-regularized functional proposed earlier. To achieve this, an accurate solution of the corresponding direct problem is essential. We tackle this challenge through a finite-element strategy for spatial discretization and a second-order implicit finite difference scheme for time discretization. The implementation is facilitated by the Dolfin-FeniCS libraries \cite{Logg, Fenics2, Farrell, Funke, Alnaes1, Alnaes1b, Alnaes2, Logg3, Fenics}, while the least-squares minimization of the Tikhonov functional is realized through the iterative L-BFGS-B routine \cite{Byrd, Nocedal, Zhu}, as implemented in the Dolfin-Adjoint libraries \cite{Farrell, Funke} and SciPy \cite{scipy} on Python.

Previous works on inverse problems related to dispersive-type equations include studies such as \cite{Cerpa2013}, \cite{Sakthivel2016}, and \cite{Sakthivel2017}, which focus specifically on the problem of identifying a coefficient in the Korteweg-de Vries equation. To the best of the authors' knowledge, this paper is the first to conduct an analytical and numerical investigation into an inverse problem for a system of dispersive equations like \eqref{KdVe1}. The primary challenge in establishing theoretical results in this context lies in the extensive calculations required to derive the Fr\'echet derivative of the Tikhonov functional, as well as the formulation of a suitable nonlinear adjoint problem essential to the process.

At this point, we establish notation to be used throughout this document. Let $L$ and $\beta$ donote positive constants, and $L^2([0, L])$ represent the space of measurable functions $u : [0, L] \to \mathbb{R}$ that are square-integrable over $[0, L]$, with its usual inner product and norm. For the Sobolev space $H^1([0,L])$, we define the following weighted inner product and induced norm:
\begin{align}
&\left\langle u,v \right\rangle_{H^1} = \int_0^L u(x)v(x) dx + \frac{\beta}{6}  \int_0^L u'(x)v'(x) dx \label{eq:innerproduct} \\
&\left\| u\right\|_{H^1} = \Big(  \int_0^L (u(x))^2 dx + \frac{\beta}{6}  \int_0^L (u'(x))^2 dx \Big)^{1/2}. \label{eq:H1norm}
\end{align}

\noindent For $f \in H^1([0,L])$, we have the inequality
\begin{equation}
\|f'\|_{L^2} \leq \frac{1}{ \Big( \frac{\beta}{6} \Big)^{\frac12} } \|f\|_{H^1}, 
\end{equation}
and for $f \in H_0^1([0,L])$ (the space of functions $f \in H^1([0,L])$ such that
$f(0) = f(L) = 0$ ), we also have
\begin{equation}
\| f \|_{L^\infty } \leq \frac{L^{\frac12} }{ \Big( \frac{\beta}{6} \Big)^{\frac12} } \|f\|_{H^1}.
\end{equation}

\noindent The following symbols will be used to denote various spaces:
\begin{equation}\label{notacion1}
	\mathcal{L}^2_{T} = L^2(0,T;H_0^1([0,L])), \quad \mathcal{C}_{T} = C(0,T;H_0^1([0,L])) \quad \text{and} \quad L^2_{T} = L^2(0,T;L^2([0,L])).
\end{equation}
In this framework, the space $\mathcal{L}^2_{T}$ is a Banach space equipped with the norm
\begin{equation*}
	\| u \|_{\mathcal{L}_T^2} = \Big( \int_0^T \| u(t) \|_{H_0^1}^2 dt \Big)^{1/2},
\end{equation*}
where the norm on $H_0^1$ is defined in \eqref{eq:H1norm}. Similarly, the spaces $\mathcal{C}_T$ and $L^2_T$ are Banach spaces with the norms given by
\begin{equation*}
	\| u \|_{\mathcal{C}_T} = \max_{t \in [0,T]} \| u(t) \|_{H_0^1} \quad \text{and} \quad \| u \|_{L^2_T} = \left( \int_0^T \| u(t) \|_{L^2}^2 dt \right)^{1/2},
\end{equation*}
respectively.

This manuscript is organized as follows.  In section 2, we formulate the direct problem considered for system \eqref{KdVe1} and establish its well-posedness and some important inequalities.  In section 3, the inverse problem is formulated and a local stability result is derived by using the Fr\'echet derivative of an appropriate regularized Tikhonov functional. In section 4, we present some numerical simulations that illustrate our analytical results, and finally, section 5 contains the conclusions of the work.

\section{Direct Problem}\label{capprodirecto}

To analyze the inverse problem effectively, it is essential to first examine the associated direct problem. In this section, we consider a more general Boussinesq system than \eqref{KdVe1}, which encompasses certain auxiliary problems necessary for the analysis.

Let $L$, $\beta$, $T$ be positive constants, and consider the system
\begin{equation}\label{Dee}
\begin{cases}
\begin{array}{c l}
N_t+\alpha\left(\xi\right)  V_{\xi}-\frac{\beta}{6}N_{\xi\xi t}&=0 \:\:\:,\:\:\left(\xi,t\right)\in \left[0,L\right]\times\left[0,T\right].\\
V_t+(c(\xi) N)_{\xi}-\frac{\beta}{6}V_{\xi\xi t}&=f(\xi,t),
\end{array}
\end{cases}
\end{equation}
subject to the initial conditions
\begin{equation}\label{Maa}
\begin{array}{c l}
N(\xi,0)&=N_0(\xi),\\
V(\xi,0)&=V_0(\xi),\:\:\: \xi \in \left[0,L\right],
\end{array}
\end{equation}
and the Dirichlet boundary conditions
\begin{equation}\label{BBMLin33}
\begin{array}{c l}
N(0,t)=N(L,t)=0,\:\: t\in \left[0,T\right]\\
V(0,t)=V(L,t)=0,\:\: t\in \left[0,T\right].
\end{array}
\end{equation}
We further assume that $c(\xi), \alpha(\xi) \in L^{2}([0,L])$ and $f(\xi,t)\in \mathcal{L}^{2}_{T}$.\\

To establish the existence of solutions for system  $\left(\ref{Dee}\right)$, we employ the Green's function associated with the differential operator $I - \frac{\beta}{6} \partial^2_{\xi}$. This allows us to reformulate problem $\left(\ref{Dee}\right)$ as a fixed point problem, to which we apply the Banach's fixed point theorem. 

The Green's function $G(\xi,s)$ associated with the differential operator $P=I-\frac{\beta}{6}\partial_{\xi}^2 $ is given by 
\begin{equation}\label{FuncionG}
G(\xi, s) = \frac{1}{2 \sqrt{\beta/6}} \frac{ \cosh\Big( \frac{1}{\sqrt{\beta/6}} (L - |s - \xi | ) \Big) - \cosh \Big( \frac{1}{\sqrt{\beta/6}} (L - (\xi + s) ) \Big)}{\sinh\Big( \frac{L}{\sqrt{\beta/6} } \Big)}.
\end{equation}
Furthermore, we define
\begin{equation}\label{DerivadaG}
K(\xi,s) :=  \frac{\partial G(\xi, s)}{\partial s} = \frac{3}{\beta} \left( \frac{ \sinh \Big( \frac{L-\xi-s}{\sqrt{\beta/6}} \Big) + \text{sign} (\xi- s) \sinh \Big( \frac{L- |\xi-s|}{\sqrt{\beta/6}} \Big)}{ \sinh \Big( \frac{L}{\sqrt{\beta/6}} \Big) }  \right), ~~ \xi \neq s.
\end{equation}
Note that
\[
\lim_{s \to \xi^-} K(\xi, s) - \lim_{s \to \xi^+} K(\xi,s) = \frac{6}{\beta}.
\]
Thus a solution $\left(N,V \right) $ of the system $\left(\ref{Dee}\right)$ must satisfy
\begin{eqnarray}\label{Ma11}
N_t & = &-\int_{0}^{L}G(\xi,s) \alpha\left(s\right) \left[  V\left( s,t\right) \right]_s ds \\
& = & \int_{0}^{L}\left( G(\xi,s)\alpha \left(s \right) \right)_{s}   V\left( s,t\right) ds, \label{Ma11new}
\end{eqnarray}

\begin{eqnarray}\label{Ma1}
V_t & = &-\int_{0}^{L}G(\xi,s) \left[ \left( c\left(s \right) N\left( s,t\right)\right)_s -f(s,t) \right] ds \\
& = & \int_{0}^{L}G_s(\xi,s) c\left( s\right) N\left( s,t\right) ds + \int_{0}^{L}G(\xi,s)f(s,t) ds.  \label{Ma1new}
\end{eqnarray}

The following lemma provides key properties of the functions $G(\xi, s)$ y $K(\xi,s)$.

\begin{lema}\label{MA51}
Let $\Phi_1$ and $\Phi_2$ be operators defined by
\begin{equation*}
\Phi_1(\phi)(\xi)=\int_{0}^{L}G(\xi,s)\phi(s)ds \quad \text{and} \quad \Phi_2(\phi)(\xi)=\int_{0}^{L}K(\xi,s)\phi(s)ds.
\end{equation*}
Then $\Phi_1$ and $\Phi_2$ are continuous linear operators from $L^{2}$ to $H^{1}$. Specifically, there exist positive constants $C_1$ and $C_2$ such that
\begin{equation}\label{MA50}
\left\| \Phi_1(\phi)\right\|_{H^{1}}\leq C_1\left\| \phi\right\|_{L^{2}}  \quad \text{and} \quad \left\| \Phi_2(\phi)\right\|_{H^{1}}\leq C_2\left\| \phi\right\|_{L^{2}}, 
\end{equation}
for all $\phi \in L^{2}$.
\end{lema}

The following theorem ensures the local well-posedness of the Boussinesq system $\left(\ref{Dee}\right)$, which can be demonstrated using the integral formulation given in \eqref{Ma11}, \eqref{Ma1}.

\begin{teorema}\label{MA3}
Let $L$, $T$ be positive constants and $f \in \mathcal{L}_T^2$. Let us consider $N_0(\xi)$, $N_0(\xi)$ $\in H_{0}^1$ and $\alpha \in  H^1$, $c \in L^2\left[0,L \right]$. Then there exists a solution $(N,V) \in \mathcal{L}_{T}^2 \times \mathcal{L}_{T}^2$ satisying $\left(\ref{Dee}\right)$, $(\ref{Maa})$, $(\ref{BBMLin33})$. This solution is unique, depends continuously on the initial data $(N_0,V_0)$ and satisfies the estimate
\begin{equation}\label{D1}
\left\| (N,V)\right\|_{\mathcal{L}_{T}^2\times \mathcal{L}_{T}^2}\leq \left(\left\| (N_0,V_0)\right\|_{H_{0}^1\times H_{0}^1} +T^{\frac{1}{2}}\left\|f\right\|_{L_{T}^2}\right)T^{\frac{1}{2}} \text{exp} \left[\frac{L^{\frac{1}{2}}}{\left( \frac{\beta}{6}\right) }D ~T \right],
\end{equation}
where $D>0$ is a constant.
\end{teorema}

\section{Inverse Problem}

In this chapter, we analyze the inverse problem of recovering the linear velocity coefficient in the Boussinesq system \eqref{Dee}-\eqref{BBMLin33} from a measure of the solution at a given final time.

\subsection{Energy estimates for the linear Boussinesq system}

Let $L$, $\beta $ be positive constants. We consider the Boussinesq system given by

\begin{equation}\label{De}
\begin{cases}
\begin{array}{c l}
N_t+V_{\xi}-\frac{\beta}{6}N_{\xi\xi t}&=0 \:\:\: ,\:\:\left(\xi,t\right)\in \left[0,L\right]\times\left[0,T\right]\\
V_t+(c(\xi) N)_{\xi}-\frac{\beta}{6}V_{\xi\xi t}&=0,
\end{array}
\end{cases}
\end{equation}
subject to initial conditions
\begin{equation}\label{Ma}
\begin{array}{c l}
N(\xi,0)&=N_0(\xi),\\
V(\xi,0)&=V_0(\xi),\:\:\: \xi \in \left[0,L\right],
\end{array}
\end{equation}
and Dirichlet boundary conditions
\begin{equation}\label{BBMLin3}
\begin{array}{c l}
N(0,t)=N(L,t)=0,\:\: t\in \left[0,T\right]\\
V(0,t)=V(L,t)=0,\:\: t\in \left[0,T\right].
\end{array}
\end{equation}

\begin{propo} (Energy estimate)\label{FE3}
Suppose that $c \in L^2\left(0,L \right) $ and $\left(N_0, V_0 \right) \in H_{0}^{1}\times H_{0}^{1}$. Then there exists a positive constant $K_c$ such as the solution to $\left(\ref{De}\right),\left(\ref{Ma}\right),\left(\ref{BBMLin3}\right)$ satisfies

\begin{equation}\label{Ana}
\left\|\left(N, V \right)  \right\|^2_{ \mathcal{L}_T^2 \times \mathcal{L}_T^2 } \leq T^{\frac{1}{2}} \left\|\left(N_0, V_0 \right)  \right\|^2_{H^{1}\times H^{1}} \exp\left( K_c T\right),
\end{equation}
where

\begin{equation}\label{MA7}
K_c=\frac{L^{\frac{1}{2}}}{\left(\frac{\beta}{6} \right) }\left\| c \right\|_{L^{2}}+\frac{1}{\left(\frac{\beta}{6} \right)^{\frac{1}{2}} }.
\end{equation}
\end{propo}

\subsection{Restricted minimization problem (RMP) }\label{sec:SecMinimization}

In this section, we reformulate the inverse problem considered as a restricted minimization problem for an appropriate Tikhonov regularized functional. It is established in first place that the minimization problem studied has at least one solution in an admisible coefficient space $\mathcal{M}_\gamma$. Consequently, for a given time $T>0$, it is possible to reformulate the inverse problem in terms of the input-output functional 
\[
\phi(c) = (N(\xi, T), V(\xi,T)),
\]
where $(N(\xi,T), V(\xi,T))$ is the solution of problem \eqref{De}-\eqref{BBMLin3} evaluated at the final time $t = T$.
This operator is well-defined by virtue of Theorem \ref{MA3}. In this way, we can rewrite the inverse problem as one of finding the coefficient $c$ such as
\begin{equation}\label{MA26}
\phi(c)=m,
\end{equation}
given a measurement $m = (m_1, m_2) \in H^{1}_{0}(\left[0,L \right] )\times  H^{1}_{0}(\left[0,L \right] )$
of the solution $(N,V)$ at the final time $t=T$. 

The set
\begin{equation*}
\mathcal{M}_\gamma=\left\lbrace c\in L^{2}(0,L): \left\| c\right\|_{L^{2}}\leq \gamma \right\rbrace, 
\end{equation*}
is the class of admissible coefficients. Now given a measurement $m=(m_1,m_2)\in  H^{1}_{0}(\left[0,L \right] )\times  H^{1}_{0}(\left[0,L \right] )$ and $\alpha \geq 0$, we define the functional $J_\alpha:\mathcal{M}_\gamma \rightarrow \mathbb{R} $ as
 
\begin{equation}\label{De7}
J_\alpha(c)= \frac{1}{2}\left\|N(\cdot, T; c) -m_1(\cdot)\right\|_{H^1}^2+ \frac{1}{2}\left\|V(\cdot, T; c) -m_2(\cdot)\right\|_{H^1}^{2}+\frac{\alpha}{2}\left\| c(\cdot)\right\| _{L^2}^{2},
\end{equation}
where the norm in the Sobolev space $H^{1}(\left[0,L \right] )$ is given in (\ref{eq:H1norm}). If $\alpha >0$ (regularization parameter), $J_\alpha$ is called a Tikhonov regularized functional.

The \textbf{Restricted Minimization Problem (RMP)} associated is the problem of minimizing the Tikhonov functional \eqref{De7}, where  $m=(m_1,m_2)\in  H^{1}_{0}(\left[0,L \right] )\times  H^{1}_{0}(\left[0,L \right] )$, $\alpha \geq 0$ and $c \in \mathcal{M}_\gamma $. In other words, we consider the problem
$$ \inf_{c\in \mathcal{M}_{\gamma}}J_\alpha (c),$$
subject to $(N(\xi,t;c),V(\xi,t;c))\in H^{1}_{0}(\left[0,L \right] )\times  H^{1}_{0}(\left[0,L \right] )$, and $N(\xi,0)=N_{0}(\xi) , V(\xi,0)=V_{0}(\xi)$, for all $\xi \in \left[ 0,L\right]$, such that the integral equations
\begin{equation*}
	\label{Ana1}
	\int_{0}^{L}N_t \phi d\xi -\int_{0}^{L}V \phi_{\xi} d\xi +\frac{\beta}{6}\int_{0}^{L}N_{\xi t} \phi_{\xi} d\xi= 0
\end{equation*}
\begin{equation*}
	\label{Ana2}
	\int_{0}^{L}V_t \psi d\xi -\int_{0}^{L}(c N)\psi_{\xi} d\xi +\frac{\beta}{6}\int_{0}^{L}V_{\xi t} \psi_{\xi} d\xi= 0,
\end{equation*}
satisfy for all $t \in [0,T]$ and $\phi,\psi \in H^{1}_{0}(\left[0,L \right] $.\\
\begin{teorema}
Let $m_1,m_2 \in H_0^1$. Then, there exists a minimizer $c^* \in \mathcal{M}_{\gamma}$ of the functional $J_\alpha$,
such that
\begin{equation*}
J_\alpha(c^*)= \inf_{c\in \mathcal{M}_{\gamma}} J_\alpha (c).
\end{equation*}
\end{teorema}

\textit{Proof:} Given that $J_\alpha(c) \geq 0$ for all $c \in \mathcal{M}_{\gamma}$, then the functional $J_\alpha(c)$ is bounded bellow. Let $\left\lbrace c_k\right\rbrace _k \subseteq \mathcal{M}_{\gamma} $ be a minimizing sequence of $J_\alpha$ such that

\begin{equation}\label{Ana5}
\lim_{k \rightarrow \infty}J_\alpha (c_k)= \inf_{c\in \mathcal{M}_{\gamma}}J_\alpha (c).
\end{equation}
Noting that $\left\|c_k \right\|_{L^2}\leq \gamma $, for all $k\in \mathbb{N}$, there exists a subsequence
in $\left\lbrace c_k\right\rbrace_k $, denoted by $\left\lbrace c_k\right\rbrace_k \subseteq \mathcal{M}_{\gamma}$, and $c^* \in L^2$, such that $c_k \rightharpoonup c^*$ in $L^2$. Observe that $\mathcal{M}_{\gamma}$ is closed in the
strong topology and convex in $L^2$, from which it is weakly-closed, i.e.  $c^* \in \mathcal{M}_{\gamma}$.\\

\noindent From $\left(\ref{Ana}\right)$, the solution $\left(N(\xi,t;c_k),V(\xi,t;c_k) \right)=\left(N(c_k),V(c_k) \right)  $ of $\left(\ref{De}\right)$ satisfies $$\left\|(N(c_k),V(c_k)) \right\|_{L^2(H^1)\times L^2(H^1)}\leq C,$$ for all $k\in \mathbb{N}$. Therefore, there exits a subsequence of $\left\lbrace \left(N(c_k),V(c_k) \right)\right\rbrace_k $, also denoted $\left\lbrace \left(N(c_k),V(c_k) \right)\right\rbrace_k $, such that 
$$N(\xi,t;c_k) \rightharpoonup N^*(\xi,t), V(\xi,t;c_k) \rightharpoonup V^*(\xi,t),$$ 
in $L^2(0,T; H^1)$. Furthermore, by Theorem \ref{MA3}, and Lions-Aubin theorem, 
we obtain strong convergence of the subsequence in $L^2(0,T; L^2(0,L))$, i.e.,

$$N(\xi,t;c_k) \rightarrow N^*(\xi,t)$$
$$V(\xi,t;c_k) \rightarrow V^*(\xi,t).$$
Next, we will show that $N^*(\xi,t)=N(\xi,t; c^*)$, $V^*(\xi,t)=V(\xi,t; c^*)$
satisfy $\left(\ref{De}\right)$ with the initial conditions $\left(\ref{BBMLin3}\right)$. 
We say that the pair $(N,V)$ is a solution to $\left(\ref{De}\right)$ if

\begin{equation}
\label{Ana1}
\int_{0}^{L}N_t \phi d\xi -\int_{0}^{L}V \phi_{\xi} d\xi +\frac{\beta}{6}\int_{0}^{L}N_{\xi t} \phi_{\xi} d\xi= 0
\end{equation}
\begin{equation}
\label{Ana2}
\int_{0}^{L}V_t \phi d\xi -\int_{0}^{L}(c N)\phi_{\xi} d\xi +\frac{\beta}{6}\int_{0}^{L}V_{\xi t} \phi_{\xi} d\xi= 0,
\end{equation}
for all $\phi \in H_0^1$ and
\begin{equation}
\begin{array}{c l}
N(\xi,0)=N_0(\xi), \\\\
V(\xi,0)=V_0(\xi).
\end{array}
\end{equation}
In fact, let $h_1 \in \tilde{C}^1=\left\lbrace h_1 \in C^{1}(\left[0, T\right]): h_1(T)=0\right\rbrace $. 
By multiplying the equation $\left(\ref{Ana1}\right)$ by $h_1$ and integrating over the interval $\left(0,T\right) $, 
we obtain
\begin{equation}
\label{Ana3}
\int_{0}^{T}\int_{0}^{L}N_t(c_k) \phi h_1 d\xi dt -\int_{0}^{T}\int_{0}^{L}V(c_k) \phi_{\xi}h_1 d\xi dt +\frac{\beta}{6}\int_{0}^{T}\int_{0}^{L}N_{\xi t}(c_k) \phi_{\xi}h_1 d\xi dt= 0.
\end{equation}
By integration by parts and rearranging terms, we obtain
\begin{equation*}
\begin{array}{c l}
-\int_{0}^{L}N_0(\xi) \phi h_1(0) d\xi - \frac{\beta}{6}\int_{0}^{L}N_{0}^{'}(\xi) \phi_{\xi}h_1(0) d\xi- \int_{0}^{T}\int_{0}^{L}V(c_k) \phi_{\xi}h_1 d\xi dt\\\\
-\int_{0}^{T}\int_{0}^{L}N(c_k) \phi h_1^{'} d\xi dt- \frac{\beta}{6}\int_{0}^{T}\int_{0}^{L}N_{\xi}(c_k) \phi_{\xi}h_1^{'} d\xi dt=0.
\end{array}
\end{equation*}
Thus
\begin{equation*}
\begin{array}{c l}
h_1(0)\left( \int_{0}^{L}N_0(\xi) \phi  d\xi + \frac{\beta}{6}\int_{0}^{L}N_{0}^{'}(\xi) \phi_{\xi} d\xi\right) =\\\\ -\int_{0}^{T}\int_{0}^{L}V(c_k) \phi_{\xi}h_1 d\xi dt -
\int_{0}^{T}\int_{0}^{L}N(c_k) \phi h_1^{'} d\xi dt- \frac{\beta}{6}\int_{0}^{T}\int_{0}^{L}N_{\xi}(c_k) \phi_{\xi}h_1^{'} d\xi dt.
\end{array}
\end{equation*}
Using that $N(\xi,t;c_k) \rightharpoonup N^*(\xi,t)$, $V(\xi,t;c_k) \rightharpoonup V^*(\xi,t)$ in $L^{2}(0,T; H^1)$, 
we obtain
\begin{equation*}
\begin{array}{c l}
h_1(0)\left( \int_{0}^{L}N_0(\xi) \phi  d\xi + \frac{\beta}{6}\int_{0}^{L}N_{0}^{'}(\xi) \phi_{\xi} d\xi\right) =\\\\ -\int_{0}^{T}\int_{0}^{L}V^{*} \phi_{\xi}h_1 d\xi dt-
\int_{0}^{T}\int_{0}^{L}N^{*} \phi h_1^{'} d\xi dt- \frac{\beta}{6}\int_{0}^{T}\int_{0}^{L}N^{*}_{\xi} \phi_{\xi}h_1^{'} d\xi dt.
\end{array}
\end{equation*}
By applying integration by parts to the last two terms, we obtain

\begin{equation*}
\begin{array}{c l}
h_1(0) \left[ \int_{0}^{L}\left( N_0(\xi)-N_{0}^{*}(\xi) \right) \phi  d_\xi - \frac{\beta}{6}\int_{0}^{L}\left( N_{0}^{'}(\xi)-N_{\xi}^{*}(\xi,0)\right)  \phi_{\xi} d_\xi\right]  =\\\\ \int_{0}^{T}\left[\int_{0}^{L}N^{*}_{t} \phi  d_\xi -\int_{0}^{L}V^{*} \phi_{\xi} d_\xi + \frac{\beta}{6}\int_{0}^{L}N^{*}_{\xi t} \phi_{\xi}d_\xi\right] h_1(t) dt,
\end{array}
\end{equation*}
for some $h_1 \in \tilde{C}^{1}\left( \left[0,T \right]\right)  $ y $\phi \in H_{0}^1$. Furthermore, rearranging terms, we yield that
\begin{equation*}
\begin{array}{c l}
h_1(0) \left( N_0(\cdot)-N_{0}^{*}(\cdot), \phi(\cdot)\right)_{H^1} =\\\\ \int_{0}^{T}\left[\int_{0}^{L}N^{*}_{t}\phi  d\xi -\int_{0}^{L}V^{*} \phi_{\xi} d\xi + \frac{\beta}{6}\int_{0}^{L}N^{*}_{\xi t} \phi_{\xi}d\xi\right] h_1(t) dt.
\end{array}
\end{equation*}
Thus, selecting $h_1 \in \tilde{C}^{1}\left( \left[0,T \right]\right)$ such that $h_1(0)=0$, we get

\begin{equation*}
\begin{array}{c l}
\int_{0}^{L}N^{*}_{t}\phi  d\xi -\int_{0}^{L}V^{*} \phi_{\xi} d_\xi + \frac{\beta}{6}\int_{0}^{L}N^{*}_{\xi t} \phi_{\xi}d_\xi =0.
\end{array}
\end{equation*}
Thus, we obtain
\begin{equation*}
\begin{array}{c l}
h_1(0) \left( N_0(\cdot)-N_{0}^{*}(\cdot), \phi(\cdot)\right)_{H^1} = 0,
\end{array}
\end{equation*}
for all $h_1 \in \tilde{C}^{1}\left( \left[0,T \right]\right)$ and $\phi \in H_{0}^1$. Then, we have that
$$\left( N_0(\cdot)-N_{0}^{*}(\cdot), \phi(\cdot)\right)_{H^1}=0,$$
for any $\phi \in H_{0}^1$. Therefore,
\begin{equation*}
\begin{array}{c l}
 N_0(\xi)= N^{*}(\xi,0).
\end{array}
\end{equation*}
Applying a similar argument to the second equation \eqref{Ana2}, and considering $$h_2 \in \tilde{C}^1=\left\lbrace h_2 (t)\in C^{1}(\left[0, T\right]): h_2(T)=0\right\rbrace,$$ 
we find that
$$\left( V_0(\cdot)-V_{0}^{*}(\cdot), \phi(\cdot)\right)_{H^1}=0,$$ 
for any $\phi \in H_{0}^1$, which implies
\begin{equation*}
\begin{array}{c l}
V_0(\xi)= V^{*}(\xi,0).
\end{array}
\end{equation*}
Thus $\left( N^{*},V^{*}\right) $ is indeed a solution to $\left(\ref{De}\right)$ that satisfies the initial conditions $\left(\ref{BBMLin3}\right)$.
The last step is to demonstrate that $c^*$ is a minimizer of the functional $J_\alpha$. First, observe that
\begin{equation}
0 \leq \left\|\left(N(\cdot, T; c_k)- m_1(\cdot)\right)-\left(N(\cdot, T; c^*)- m_2(\cdot)\right) \right\|_{H^1} ^2.
\end{equation}
Therefore, given that $N(\cdot, T; c_k) \rightharpoonup N(\cdot, T; c^*)$, it follows that
\begin{equation}
\begin{array}{c l}
\left\| N(\cdot, T; c^*)- m_1(\cdot)\right\|_{H^1} ^2=2\lim_{k\rightarrow \infty} ( N(\cdot, T; c_k)- m_1(\cdot),N(\cdot, T; c^*)- m_1(\cdot) )_{H^1}\\\\
-\left\|N(\cdot, T; c^*)- m_1(\cdot) \right\|_ {H^1}^2\\\\
\leq \lim_{k\rightarrow \infty}\inf \left\|N(\cdot, T; c_k)- m_1(\cdot) \right\|_ {H^1}^2.
\end{array}
\end{equation}
For the second component $V$, we have that $V(\cdot, T; c_k) \rightharpoonup V(\cdot, T; c^*)$ in $L^2_{H^1}$, and thus
\begin{equation}
\begin{array}{c l}
\left\| V(\cdot, T; c^*)- m_2(\cdot)\right\|_{H^1} ^2=2\lim_{k\rightarrow \infty}\left(V(\cdot, T; c_k)- m_2(\cdot),V(\cdot, T; c^*)- m_2(\cdot) \right)_{H^1}\\\\
-\left\|V(\cdot, T; c^*)- m_2(\cdot) \right\|_ {H^1}^2\\\\
\leq \lim_{k\rightarrow \infty}\inf \left\|V(\cdot, T; c_k)- m_2(\cdot) \right\|_ {H^1}^2.
\end{array}
\end{equation}
Therefore, using the weak convergence $c_k \rightharpoonup c^*$ in $L^2$ and the lower semicontinuity in $H^1$, it follows that 
\begin{equation*}
\begin{array}{c l}
J_\alpha(c^*)=\frac{1}{2}\left\|N(\cdot, T; c^*) -m_1(\cdot)\right\|_{H^1}^2+ \frac{1}{2}\left\|V(\cdot, T; c^*) -m_2(\cdot)\right\|_{H^1}^{2}+\frac{\alpha}{2}\left\| c^*(\cdot)\right\| _{L^2}^{2}\\\\
\leq \frac{1}{2}\lim_{k\rightarrow \infty}\inf \left\|N(\cdot, T; c_k)- m_1(\cdot) \right\|_ {H^1}^2+\frac{1}{2}\lim_{k\rightarrow \infty}\inf \left\|V(\cdot, T; c_k)- m_2(\cdot) \right\|_ {H^1}^2+ \frac{\alpha}{2}\left\| c_k(\cdot)\right\| _{L^2}^{2}\\\\
= \lim_{k\rightarrow \infty}\inf J_\alpha (c_k).
\end{array}
\end{equation*}
Now using $\left(\ref{Ana5}\right)$, we obtain
\begin{equation*}
\begin{array}{c l}
\inf_{c\in \mathcal{M_\gamma}} J_\alpha (c)=\lim_{k \rightarrow \infty}J_\alpha (c_k)=\liminf_{k \rightarrow \infty} J_\alpha (c_k)\geq  J_\alpha (c^*)\geq \inf_{c\in \mathcal{M}_\gamma} J_\alpha(c)
\end{array}
\end{equation*}
\begin{equation*}
\begin{array}{c l}
\inf_{c\in \mathcal{M_\gamma}} J_\alpha (c)= J_\alpha (c^*).
\end{array}
\end{equation*}
Finally, we conclude that $c^*$ is a minimizer of the functional $J_\alpha$ in $\mathcal{M_{\gamma}}$.

$\hfill\square$

\subsection{Study of the gradient and optimality conditions}\label{sec:GradFrechetycondopt}

In this subsection, we derive an explicit formula for the Fréchet derivative of the regularized Tikhonov functional $\mathcal{F}_{\alpha}$ and establish the first-order optimality conditions for the solutions of the restricted minimization problem (RMP).

Let $\tilde{c}, c$ in $L^{2}(\left[ 0,L\right] )$ be two coefficients, and let $\tilde{U}=(\tilde{N},\tilde{N})$, $U=(N,V)$ denote the solutions to the problem (\ref{De}) corresponding to these coefficients $\tilde{c}$ and $c$, respectively, i.e., $\tilde{U}=U(\tilde{c})$ and $U=U(c)$. We introduce the functions $\mathcal{N}(\xi,t;\tilde{c},c)$, $\mathcal{V}(\xi,t;\tilde{c},c)$ defined by
\begin{equation}
\mathcal{N}(\xi,t;\tilde{c},c)=\mathcal{N}(\xi,t):=\tilde{N}(\xi,t)-N(\xi,t)
\end{equation}
\begin{equation}
\mathcal{V}(\xi,t;\tilde{c},c)=\mathcal{V}(\xi,t):=\tilde{V}(\xi,t)-V(\xi,t).
\end{equation}
Thus, the pair $(\mathcal{N},\mathcal{V})$ satisfies the system
\begin{equation}\label{MA1}
\begin{cases}
\begin{array}{c l}
\mathcal{N}_t+ \mathcal{V}_{\xi}-\frac{\beta}{6}\mathcal{N}_{\xi\xi t}&=0 \\
\mathcal{V}_t+(c \mathcal{N})_{\xi}-\frac{\beta}{6}\mathcal{V}_{\xi\xi t}&=-((\tilde{c}-c)\tilde{N})_{\xi},\\
\mathcal{N}(\xi,0)=\mathcal{V}(\xi,0)=0\\
\mathcal{V}(0,t)=\mathcal{V}(L,t)=0\\
\mathcal{N}(0,t)=\mathcal{N}(L,t)=0.
\end{array}
\end{cases}
\end{equation}

\begin{obs}
Observe that the initial-boundary value problem \eqref{MA1} is in the same form as the system \eqref{Dee}, 
ensuring well-posedness in the space $\mathcal{L}_{T}^{2} \times \mathcal{L}_T^2$.
\end{obs}

The following proposition provides an energy estimate for system (\ref{MA1}).
\begin{propo}\label{MA11}
Let $(N_0,V_0) \in H_0^1 \times H_0^1$, $(m_1,m_2)$ $\in H^{1}_{0} \times H_0^1$, $\tilde{c},c \in L^{2}(\left[ 0,L\right] )$. Then the unique solution $(\mathcal{N},\mathcal{V})$ to problem (\ref{MA1}) satisfies the estimate
\begin{equation}\label{MA8}
\left\|(\mathcal{N},\mathcal{V})\right\|_{\mathcal{L}_{T}^{2}\times \mathcal{L}_{T}^{2}}\leq  \frac{L^{\frac{1}{2}}T^{\frac{3}{2}}}{(\frac{\beta}{6})} \left\|\tilde{c}-c\right\|_{L^{2}}\left\|(N_0,V_0) \right\|_{H^{1}\times H^{1}}K_1(\tilde{c},c),
\end{equation} 
where
\begin{equation}
K_{1}(c,\tilde{c})=\exp \left[\left( K_{\tilde{c}}+K_{c}\right)  T\right],
\end{equation}
and
\begin{equation}
K_{c}=\frac{L^{\frac{1}{2}}}{\left( \frac{\beta}{6}\right)}\left\|c \right\|_{L^{2}}+ \frac{1}{\left( \frac{\beta}{6}\right)^{\frac{1}{2}}}.
\end{equation}

\end{propo}

\textit{Proof:} Multiplying the first equation in system (\ref{MA1}) by $\mathcal{N}$, integrating over $\left[ 0,L\right]$, 
and applying integration by parts, we obtain
\begin{equation}\label{MA2}
\int_{0}^{L}\mathcal{N}\mathcal{N}_t + \frac{\beta}{6}\int_{0}^{L}\mathcal{N}_\xi \mathcal{N}_{\xi t }d\xi =\int_{0}^{L} \mathcal{N}_\xi \mathcal{V} d\xi,
\end{equation}
where 
\begin{equation*}
\int_{0}^{L}\mathcal{N}\mathcal{N}_t+\frac{\beta}{6}\int_{0}^{L}\mathcal{N}_\xi \mathcal{N}_{\xi t}d_\xi =\frac{1}{2}\frac{d}{dt}\left\| \mathcal{N}(t)\right\|^{2}_{H_{1}}.
\end{equation*}
Substituting into (\ref{MA2}), and using similar estimates as in the proof of
Theorem (\ref{MA3}), we get
\begin{equation*}
\begin{split}
\frac{1}{2}\frac{d}{dt}\left\| \mathcal{N}(t)\right\|^{2}_{H_{1}}&\leq \left\|\mathcal{N}_\xi \right\|_{L^{2}} \left\|\mathcal{V} \right\|_{L^{2}}\\
&\leq \frac{1}{(\frac{\beta}{6})^{\frac{1}{2}}}\left\|\mathcal{N} \right\|_{H^{1}} \left\|\mathcal{V} \right\|_{H^{1}},
\end{split}
\end{equation*}
\begin{equation*}
\begin{split}
\frac{d}{dt}\left\| \mathcal{N}(t)\right\|_{H_{1}}&\leq \frac{1}{(\frac{\beta}{6})^{\frac{1}{2}}} \left\|\mathcal{V} \right\|_{H^{1}} .
\end{split}
\end{equation*}
Integrating over the interval $(0,t)$, we obtain
\begin{equation}\label{MA6}
\begin{split}
\left\| \mathcal{N}(t)\right\|_{H_{1}}&\leq  \left\|\mathcal{N}_0\right\|_{H^{1}}+ \frac{1}{(\frac{\beta}{6})^{\frac{1}{2}}} \int_{0}^{t}\left\|\mathcal{V}(\tau) \right\|_{H^{1}}d\tau.
\end{split}
\end{equation}
By following similar steps for the second equation in \eqref{MA1}, and applying Gronwall's inequality,
we ultimate arrive at
\begin{equation}
\left\|(\mathcal{N},\mathcal{V})\right\|_{\mathcal{L}_{T}^{2}\times \mathcal{L}_{T}^{2}}\leq  \frac{L^{\frac{1}{2}}T^{\frac{3}{2}}}{(\frac{\beta}{6})} \left\|\tilde{c}-c\right\|_{L^{2}}\left\|(N_0,V_0) \right\|_{H^{1}\times H^{1}} \exp\left[\left( K_{\tilde{c}}+K_{c}\right)  T\right]. 
\end{equation} 
$\hfill\square$

Next, we introduce an auxiliary system to problem (\ref{De}) associated to the coefficient $c \in L^{2}(\left[0,L \right] )$, called adjoint problem for system (\ref{De}):
\begin{equation}\label{MA9}
\begin{cases}
\begin{array}{c l}
\eta_t+ c(\xi) \gamma_{\xi}-\frac{\beta}{6}\eta_{\xi\xi t}&=0 \\
\gamma_t+\eta_{\xi}-\frac{\beta}{6}\gamma_{\xi\xi t}&=0,\\
\eta(\xi,T)=\eta_T(\xi), \gamma(\xi,T)=\gamma_T(\xi) \\
\eta(0,t)=\eta(L,t)=0\\
\gamma(0,t)=\gamma(L,t)=0,
\end{array}
\end{cases}
\end{equation}
where $\eta_T , \gamma_T \in H^{1}_{0}(\left[ 0,L\right] )$. The following proposition provides an energy estimate for the adjoint system (\ref{MA9}).\\
\begin{propo}
Let $c \in L^{2}(\left[0,L \right] )$, $U=(N,V)$ be the solution to the direct problem (\ref{De})
associated to the coefficient $c$, and $\eta_T, \gamma_T \in H^{1}_{0}(\left[ 0,L\right] )$. Then the solution $(\eta,\gamma)$ to the adjoint problem (\ref{MA9}), satisfies the following energy estimate:

\begin{equation}\label{MA10}
\left\|(\eta, \gamma) \right\|_{\mathcal{L}_{T}^{2} \times \mathcal{L}_{T}^{2} }\leq T^{\frac{1}{2}}\left\| (\eta_T,\gamma_T)\right\|_{H^{1}\times H^{1}}K_2(c),
\end{equation}
where 
\begin{equation}\label{FE1}
K_2(c)= \exp (K_cT),
\end{equation}
and $K_c$ is defined as in (\ref{MA7}).
\end{propo}
\noindent This result can be proved using the same strategy as in the proof of Proposition \ref{MA11}.

\bigskip

In the next result, we derive an expression for the Fr\'echet derivative of the regularized Tikhonov functional $J_\alpha$ defined in $\left( \ref{De7}\right) $. \\

\begin{teorema}
Let $(N_0,V_0)$, $(m_1,m_2)$ $\in H_{0}^{1} \times H_0^1$ and let $c,\tilde{c}$ be admissible coefficients in
 $\mathcal{M}_{\gamma}$. Furthermore, let $\mathcal{N}=N(\tilde{c})-N(c)=\tilde{N}-N$, $\mathcal{V}=V(\tilde{c})-V(c)=\tilde{V}-V$ and $(\eta,\gamma)$ be solutions to problems (\ref{MA1}), (\ref{MA9}), respectively. 
 Then the regularized Tikhonov functional $J_\alpha$ is Fr\'echet differentiable and its derivative is given by
\begin{equation}
J_{\alpha}^{'}(\xi)=\int_{0}^{T}N(\xi,t;c)\gamma_{\xi}(\xi,t;c)dt +\alpha c(\xi),
\end{equation}
where $(\eta, \gamma)$ satisfies the final data $\eta(\xi,T)=N_{T}(\xi)=N(\xi,T;c)-m_1(\xi)$, $\gamma(\xi,T)=V_{T}(\xi)=V(\xi,T;c)-m_2(\xi)$, respectively.
\end{teorema}

\textit{Proof:}
Let $c,\tilde{c} \in \mathcal{M}_\gamma$ and let $\delta J_\alpha(c)$ denote the variation of the Tikhonov
functional, i.e., $\delta J_\alpha(c)=J_\alpha(\tilde{c})-J_\alpha(c)$.  Observe that
\[
\begin{split}
J_{\alpha}(c)&=\frac{1}{2}\Vert N(\cdot,T,c)-m_1(\cdot)\Vert^{2}_{H^1_{0}}+\frac{1}{2}\Vert V(\cdot,T,c)-m_2(\cdot)\Vert^{2}_{H^1_{0}}+\frac{\alpha}{2}\Vert c \Vert^{2}_{L^2}, \\
J_{\alpha}(\tilde{c})&=\frac{1}{2}\Vert N(\cdot,T,\tilde{c})-m_1(\cdot)\Vert^{2}_{H^1_{0}}+\frac{1}{2}\Vert V(\cdot,T,\tilde{c})-m_2(\cdot)\Vert^{2}_{H^1_{0}}+\frac{\alpha}{2}\Vert\tilde{c} \Vert^{2}_{L^2}.
\end{split}
\]
Thus, we can express the variation as
\[
\begin{split}
J_\alpha(\tilde{c})-J_\alpha(c) &=\frac{1}{2}\Vert \tilde{N}(T)-m_1\Vert^2_{H^1_{0}}+\frac{1}{2}\Vert \tilde{V}(T)-m_2\Vert^2_{H^1_{0}}+\frac{\alpha}{2}\Vert\tilde{c} \Vert ^2_{L^2} \\
&\ -\frac{1}{2}\Vert N(T)-m_1\Vert ^2_{H^1_{0}}-\frac{1}{2}\Vert V(T)-m_2\Vert^2_{H^1_{0}}-\frac{\alpha}{2}\Vert c \Vert^2_{L^2}.
\end{split}
\]
To estimate the right-hand side of the previous equation, note that
\begin{equation}
\begin{split}
\frac{\alpha}{2}\Vert \tilde{c} \Vert^2_{L^2}-\frac{\alpha}{2}\Vert c \Vert^2_{L^2}&=
\frac{\alpha}{2} \left[\left\langle \tilde{c},\tilde{c}\right\rangle _{L^2}-\left\langle c,c\right\rangle _{L^2}\right]\\
&= \frac{\alpha}{2} \left[\left\langle c,c\right\rangle _{L^2}+\left\langle c,\tilde{c}-c\right\rangle _{L^2}+\left\langle \tilde{c}-c,c\right\rangle _{L^2}+\left\langle \tilde{c}-c,\tilde{c}-c\right\rangle _{L^2}-\left\langle c,c\right\rangle _{L^2}\right]\\
 &= \frac{\alpha}{2}\left[2\left\langle c,\tilde{c}-c\right\rangle _{L^2}+\left\langle \tilde{c}-c,\tilde{c}-c\right\rangle_{L^2}\right]\\
&=\alpha\left\langle c,\tilde{c}-c\right\rangle _{L^2}+\frac{\alpha}{2}\left\langle \tilde{c}-c,\tilde{c}-c\right\rangle_{L^2}.
\end{split}
\end{equation}
Adding and subtracting convenient terms and using the substitutions $\mathcal{N}=\tilde{N}-N$, $\mathcal{V}=\tilde{V}-V$ , it follows that
\begin{multline}\label{DE6}
\delta J_{\alpha}(c)=\left\langle \mathcal{N}(\cdot,T),N(\cdot,T)-m_1(\cdot)\right\rangle_{H^{1}}+\left\langle \mathcal{V}(\cdot,T),V(\cdot,T)-m_2(\cdot)\right\rangle_{H^{1}}+\alpha \left\langle c,\tilde{c}-c\right\rangle_{L^{2}}\\
+\frac{1}{2}\left\langle \mathcal{N}(\cdot,T),\mathcal{N}(\cdot,T)\right\rangle_{H^{1}}  +\frac{1}{2}\left\langle \mathcal{V}(\cdot,T),\mathcal{V}(\cdot,T)\right\rangle_{H^{1}} +\frac{\alpha}{2}\left\langle \tilde{c}-c,\tilde{c}-c\right\rangle_{L^{2}},
\end{multline}
where $(\mathcal{N},\mathcal{V})$ is solution to equations (\ref{MA1} ).

Next, let us examine the term $\left\langle \mathcal{N}(\cdot,T),N(\cdot,T)-m_1(\cdot)\right\rangle_{H^{1}}$:
\begin{equation*}
\begin{split}
\left\langle \mathcal{N}(\cdot,T),N(\cdot,T)-m_1(\cdot)\right\rangle_{H^{1}}&=\left\langle \mathcal{N}(\cdot,T),\eta_{T}(\cdot)\right\rangle_{H^{1}}\\
&=\int_{0}^{L}\mathcal{N}(\xi,T)\eta_{T}(\xi)d\xi+\frac{\beta}{6}\int_{0}^{L}\mathcal{N}_{\xi}(\xi,T)\eta_{T}^{'}(\xi)d\xi\\
&=\frac{\beta}{6}\int_{0}^{L}\mathcal{N}_{\xi}(\xi,T)\eta_{T}^{'}(\xi)d\xi+\int_{0}^{L}\eta(\xi,t)\mathcal{N}(\xi,t)\left. \right|^{t=T}_{t=0}d\xi\\
&=\frac{\beta}{6}\int_{0}^{L}\mathcal{N}_{\xi}(\xi,T)\eta_{T}^{'}(\xi)d\xi+\int_{0}^{L}\int_{0}^{T}\left( \eta(\xi,t)\mathcal{N}(\xi,t)\right)_{t}dtd\xi\\
&=\frac{\beta}{6}\int_{0}^{L}\mathcal{N}_{\xi}(\xi,T)\eta_{T}^{'}(\xi)d\xi+\int_{0}^{L}\int_{0}^{T} \eta_{t}\mathcal{N}+\eta\mathcal{N}_{t} dt d\xi. \\
\end{split}
\end{equation*}
Taking into account the expressions for $\mathcal{N}_{t}$ y $\eta_{t}$
given by equations (\ref{MA1}) y (\ref{MA9}), and integrating by parts, we obtain the following:
\begin{equation*}
\int_{0}^{L}\int_{0}^{T} \eta_{t}\mathcal{N}+\eta\mathcal{N}_{t} dt d\xi=\int_{0}^{T}\int_{0}^{L}-c\gamma_{\xi}\mathcal{N}+\frac{\beta}{6}\eta_{\xi\xi t}\mathcal{N}-\eta\mathcal{V}_{\xi}+\frac{\beta}{6}\eta\mathcal{N}_{\xi\xi t}d\xi dt.
\end{equation*}
Thus, we have
\begin{equation*}
\int_{0}^{L}\int_{0}^{T} \eta_{t}\mathcal{N}+\eta\mathcal{N}_{t} dt d\xi=\int_{0}^{T}\int_{0}^{L}-c\gamma_{\xi}\mathcal{N}-\eta\mathcal{V}_{\xi}-\frac{\beta}{6}\int_{0}^{L}\mathcal{N}_{\xi}(\xi,T)\eta ^{'}_{T}(\xi)d\xi.
\end{equation*}
From this, we obtain
\begin{equation}\label{DE2}
\begin{split}
\left\langle \mathcal{N}(\cdot,T),N(\cdot,T)-m_1(\cdot)\right\rangle_{H^{1}}
&=\int_{0}^{T}\int_{0}^{L}-c\gamma_{\xi}\mathcal{N}-\eta\mathcal{V}_{\xi} d\xi dt. \\
\end{split}
\end{equation}
Next, consider the second term on the right-hand side of equation
 (\ref{DE6}). We have
\begin{equation*}
\begin{split}
\left\langle \mathcal{V}(\cdot,T),V(\cdot,T)-m_2(\cdot)\right\rangle_{H^{1}}&=\left\langle \mathcal{V}(\cdot,T),\gamma_{T}(\cdot)\right\rangle_{H^{1}}\\
&=\int_{0}^{L}\mathcal{V}(\xi,T)\gamma_{T}(\xi)d\xi+\frac{\beta}{6}\int_{0}^{L}\mathcal{V}_{\xi}(\xi,T)\gamma_{T}^{'}(\xi)d\xi\\
&=\frac{\beta}{6}\int_{0}^{L}\mathcal{V}_{\xi}(\xi,T)\gamma_{T}^{'}(\xi)d\xi+\int_{0}^{L}\gamma(\xi,t)\mathcal{V}(\xi,t)\left. \right|^{t=T}_{t=0}d\xi\\
&=\frac{\beta}{6}\int_{0}^{L}\mathcal{V}_{\xi}(\xi,T)\gamma_{T}^{'}(\xi)d\xi+\int_{0}^{L}\int_{0}^{T}\left( \gamma(\xi,t)\mathcal{V}(\xi,t)\right)_{t}dtd\xi\\
&=\frac{\beta}{6}\int_{0}^{L}\mathcal{V}_{\xi}(\xi,T)\gamma_{T}^{'}(\xi)d\xi+\int_{0}^{L}\int_{0}^{T} \gamma_{t}\mathcal{V}+\gamma\mathcal{V}_{t} dt d\xi. \\
\end{split}
\end{equation*}
Thus, we obtain
\begin{equation*}
\int_{0}^{L}\int_{0}^{T} \gamma_{t}\mathcal{V}+\gamma\mathcal{V}_{t} dt d\xi=\int_{0}^{T}\int_{0}^{L}-\eta_{\xi}\mathcal{V}+\frac{\beta}{6}\gamma_{\xi\xi t}\mathcal{V}-\gamma(c\mathcal{N})_{\xi}+\frac{\beta}{6}\gamma\mathcal{V}_{\xi\xi t}-\gamma ((\tilde{c}-c)\tilde{N})_{\xi}d\xi dt.
\end{equation*}
From this, we get
\begin{equation*}
\int_{0}^{L}\int_{0}^{T} \gamma_{t}\mathcal{V}+\gamma\mathcal{V}_{t} dt d\xi=\int_{0}^{T}\int_{0}^{L}-\eta_{\xi}\mathcal{V}-\gamma(c\mathcal{N})_{\xi}-\gamma ((\tilde{c}-c)\tilde{N})_{\xi}d\xi dt-\frac{\beta}{6}\int_{0}^{L}\mathcal{V}_{\xi}(\xi,T)\gamma_{T}^{'}(\xi)d\xi.
\end{equation*}
Therefore, we arrive at
\begin{equation}\label{DE3}
\begin{split}
\left\langle \mathcal{V}(\cdot,T),V(\cdot,T)-m_2(\cdot)\right\rangle_{H^{1}}
&=\int_{0}^{T}\int_{0}^{L}-\eta_{\xi}\mathcal{V}-\gamma(c\mathcal{N})_{\xi}-\gamma ((\tilde{c}-c)\tilde{N})_{\xi}d\xi dt. \\
\end{split}
\end{equation}
Substituting equations (\ref{DE2}) and (\ref{DE3}) into equation \eqref{DE6}, we get
\begin{multline}
\delta J_{\alpha}(c)=\int_{0}^{T}\int_{0}^{L}-c\gamma_{\xi}\mathcal{N}-\eta\mathcal{V}_{\xi}+\int_{0}^{T}\int_{0}^{L}-\eta_{\xi}\mathcal{V}-\gamma(c\mathcal{N})_{\xi}-\gamma ((\tilde{c}-c)\tilde{N})_{\xi}d\xi dt+\alpha \left\langle c,\tilde{c}-c\right\rangle_{L^{2}}\\
+\frac{1}{2}\left\langle \mathcal{N}(\cdot,T),\mathcal{N}(\cdot,T)\right\rangle_{H^{1}}  +\frac{1}{2}\left\langle \mathcal{V}(\cdot,T),\mathcal{V}(\cdot,T)\right\rangle_{H^{1}} +\frac{\alpha}{2}\left\langle \tilde{c}-c,\tilde{c}-c\right\rangle_{L^{2}},
\end{multline}
and thus
\begin{multline*}
\delta J_{\alpha}(c)=\int_{0}^{T}\int_{0}^{L}\gamma_{\xi} ((\tilde{c}-c)\tilde{N})d\xi dt+\alpha \left\langle c,\tilde{c}-c\right\rangle_{L^{2}}\\
+\frac{1}{2}\left\langle \mathcal{N}(\cdot,T),\mathcal{N}(\cdot,T)\right\rangle_{H^{1}}  +\frac{1}{2}\left\langle \mathcal{V}(\cdot,T),\mathcal{V}(\cdot,T)\right\rangle_{H^{1}} +\frac{\alpha}{2}\left\langle \tilde{c}-c,\tilde{c}-c\right\rangle_{L^{2}},
\end{multline*}
or equivalently
\begin{multline*}
\delta J_{\alpha}(c)=\int_{0}^{T}\int_{0}^{L}\gamma_{\xi} (\tilde{c}-c)(\mathcal{N}+N)d\xi dt+\alpha \left\langle c,\tilde{c}-c\right\rangle_{L^{2}}\\
+\frac{1}{2}\left\langle \mathcal{N}(\cdot,T),\mathcal{N}(\cdot,T)\right\rangle_{H^{1}}  +\frac{1}{2}\left\langle \mathcal{V}(\cdot,T),\mathcal{V}(\cdot,T)\right\rangle_{H^{1}} +\frac{\alpha}{2}\left\langle \tilde{c}-c,\tilde{c}-c\right\rangle_{L^{2}}.
\end{multline*}
Rearranging terms, we obtain
\begin{equation*}
\delta J_{\alpha}(c)=\left\langle\int_{0}^{T}N(\cdot,t) \gamma_{\xi}(\cdot,t)dt,\tilde{c}-c\right\rangle_{L^{2}}+\alpha\left\langle c,\tilde{c}-c\right\rangle_{L^{2}}+R(\tilde{c}-c,\mathcal{N},\mathcal{V}),
\end{equation*}
where
\begin{multline*}
R(\tilde{c}-c,\mathcal{N},\mathcal{V})=\int_{0}^{T}\int_{0}^{L}\gamma_{\xi} (\tilde{c}-c)\mathcal{N}d\xi dt
+\frac{1}{2}\left\langle \mathcal{N}(\cdot,T),\mathcal{N}(\cdot,T)\right\rangle_{H^{1}}\\ +\frac{1}{2}\left\langle \mathcal{V}(\cdot,T),\mathcal{V}(\cdot,T)\right\rangle_{H^{1}} +\frac{\alpha}{2}\left\langle \tilde{c}-c,\tilde{c}-c\right\rangle_{L^{2}}.
\end{multline*}
To estimate the term $R(\tilde{c}-c,\mathcal{N},\mathcal{V})$, observe that
\begin{equation*}
\begin{split}
\left|\int_{0}^{T}\int_{0}^{L}\gamma_{\xi}\mathcal{N}(\tilde{c}-c)d\xi dt\right|&\leq \left\|\tilde{c}-c \right\|_{L^{2}}\int_{0}^{T}\left\|\mathcal{N}(t)\right\|_{L^{\infty}}\left\|\gamma_{\xi}(t)\right\|_{L^{2}}dt\\
&\leq \frac{L^{\frac{1}{2}}}{\left( \frac{\beta}{6}\right) }\left\|\tilde{c}-c \right\|_{L^{2}}\int_{0}^{T}\left\|\mathcal{N}(t)\right\|_{H^{1}}\left\|\gamma(t)\right\|_{H^{1}}dt\\
&\leq \frac{L^{\frac{1}{2}}}{\left( \frac{\beta}{6}\right) }\left\|\tilde{c}-c \right\|_{L^{2}}\left\|\mathcal{N}\right\|_{\mathcal{L}^{2}_{T}}\left\|\gamma\right\|_{\mathcal{L}^{2}_{T}}.
\end{split}
\end{equation*}
Thus, we have that estimates (\ref{MA8}) and (\ref{MA10}) imply that there exists a constant $C_{1}>0$ such that
\begin{equation*}
\begin{split}
\left|\int_{0}^{T}\int_{0}^{L}\gamma_{\xi}\mathcal{N}(\tilde{c}-c)d\xi dt\right|\leq C_{1}\left\|\tilde{c}-c \right\|^{2}_{L^{2}}.
\end{split}
\end{equation*}
As a consequence, the Fr\'echet derivative of the functional $J_{\alpha}$ in the direction $\tilde{c}-c$, can be written as
\begin{equation*}
 J_{\alpha}^{'}(c)(\tilde{c}-c)=\left\langle\int_{0}^{T}N(\cdot,t) \gamma_{\xi}(\cdot,t)dt,(\tilde{c}-c)(\cdot)\right\rangle_{L^{2}}+\alpha\left\langle c,\tilde{c}-c\right\rangle_{L^{2}}.
\end{equation*}
Therefore,
\begin{equation}\label{DE4}
J_{\alpha}^{'}(c)(\xi)=\int_{0}^{T}N(\xi,t;c) \gamma_{\xi}(\xi,t;c)dt +\alpha c(\xi).
\end{equation}
$\hfill\square$

With the help of the Fr\'echet derivative of the functional $J_{\alpha}$, the following theorem allows to establish the optimality conditions for the RMP formulation.

\begin{teorema}
Let $c\in \mathcal{M}_\gamma$ be a solution to the restricted minimization problem (RMP) and
$U=(N,V)$ be the solution to the system (\ref{De}) associated to the coefficient $c$, i.e., $U=U(c)$. Then the
following inequality satisfies for all $h\in \mathcal{M}_\gamma$,
\begin{equation}\label{MA20}
\left\langle \int^{L}_{0}N(\cdot,t)\gamma_\xi (\cdot,t)dt,(h-c)(\cdot)\right\rangle_{L^{2}}  +\alpha\left\langle c,(h-c)\right\rangle_{L^{2}}\geq 0.
\end{equation}
Here $(\eta, \gamma)$ is solution to the adjoint problem (\ref{MA9}) associated to the coefficient $c$ and with final data
\begin{equation*}
\begin{split}
\eta_T(\xi)&=N(\xi,T;c)-m_1(\xi)\\
\gamma_T(\xi)&=V(\xi,T;c)-m_2(\xi).
\end{split}
\end{equation*}

\end{teorema}
\textit{Proof:}
Let $c$ be a solution to RMP and let $U=U(c)$. Consider the following perturbation of $c$:
\begin{equation*}
c_{s}=c+s(h-c),
\end{equation*}
where $s \in \left[0,1\right] $, $h \in \mathcal{M}_\gamma$ and $\tilde{U}=(\tilde{N},\tilde{V})=(N(\xi,t;c_s),V(\xi,t;c_s))$ is the solution to the direct problem (\ref{De}) corresponding to the coefficient $c_s$. 

By considering the variation of the functional (\ref{DE6}), it follows that
\begin{multline*}
\frac{J_\alpha(c_s)-J_\alpha(c)}{s}=\left\langle\frac{\tilde{N}(\cdot,T)-N(\cdot,T)}{s},N(\cdot,T)-m_1(\cdot)\right\rangle_{H^{1}}\\+\left\langle\frac{\tilde{V}(\cdot,T)-V(\cdot,T)}{s},V(\cdot,T)-m_2(\cdot) \right\rangle_{H^{1}}\\
+\frac{\alpha}{s}\left\langle c,s(h-c)\right\rangle_{L^{2}}+\frac{1}{2s}\left\|\tilde{N}(\cdot,T)-N(\cdot,T)\right\|^{2}_{H^{1}}+\frac{1}{2s}\left\|\tilde{V}(\cdot,T)-V(\cdot,T)\right\|^{2}_{H^{1}}+\frac{\alpha}{2s}\left\|s(h-c) \right\|^{2}_{L^{2}}.
\end{multline*}
Let us define
\begin{equation*}
K_s(\xi,t)=\frac{\tilde{N}(\xi,t)-N(\xi,t)}{s}, ~~ 
F_s(\xi,t)=\frac{\tilde{V}(\xi,t)-V(\xi,t)}{s}.
\end{equation*}
Thus, the variation becomes
\begin{multline*}
\frac{J_\alpha(c_s)-J_\alpha (c)}{s}=\left\langle K_s(\cdot,T),N(\cdot,T)-m_1(\cdot)\right\rangle_{H^{1}}+\left\langle F_s(\cdot,T),V(\cdot,T)-m_2(\cdot) \right\rangle_{H^{1}}\\
+\alpha\left\langle c,(h-c)\right\rangle_{L^{2}}+\frac{s}{2}\left\|K_s(\cdot,T)\right\|^{2}_{H^{1}}+\frac{s}{2}\left\|F_s(\cdot,T)\right\|^{2}_{H^{1}}+\frac{\alpha s}{2}\left\|(h-c) \right\|^{2}_{L^{2}}.
\end{multline*}
We can establish an analogous estimate to (\ref{MA8}) using the result from proposition (\ref{MA11}). Explicitly, this yields
\begin{equation}
\left\|(K_s(\cdot,T),F_s(\cdot,T))\right\|_{H^{1}\times {H^{1}}}\leq  \frac{L^{\frac{1}{2}}T}{(\frac{\beta}{6})} \left\|h-c\right\|_{L^{2}}\left\|(N_0,V_0) \right\|_{H^{1}\times H^{1}}K_1(c_s,c).
\end{equation} 
It is important to note that the right-hand side of this equation remains bounded when $s$ tends to cero. This, together with
the fact that $J_\alpha$ is Fr\'echet differentiable at $c_s$, implies
\begin{multline}\label{MA13}
\frac{d}{ds}J_\alpha(c_s)\left|_{s=0}\right.=\left\langle K_s(\cdot,T)\left|_{s=0} \right. ,N(\cdot,T)-m_1(\cdot)\right\rangle_{H^{1}}\\+\left\langle F_s(\cdot,T)\left|_{s=0} \right.,V(\cdot,T)-m_2(\cdot) \right\rangle_{H^{1}}
+\alpha\left\langle c,(h-c)\right\rangle_{L^{2}}.
\end{multline}
A direct calculation shows that $(K_s,F_s)$ satisfies the problem
\begin{equation}
\begin{cases}
\begin{array}{c l}
(K_{s})_t+(F_{s})_{\xi}-\frac{\beta}{6}(K_{s})_{\xi\xi t}&=0 \\
(F_{s})_t+(cK_{s})_{\xi}-\frac{\beta}{6}(F_{s})_{\xi\xi t}&=((c-h)\tilde{N})_{\xi},\\
K_{s}(\xi,0)=F_{s}(\xi,0)=0\\
F_{s}(0,t)=F_{s}(L,t)=0\\
K_{s}(0,t)=K_{s}(L,t)=0.
\end{array}
\end{cases}
\end{equation}
Let us define $(\zeta,\sigma)$ as $\zeta=K_s\left|_{s=0} \right.$ and $\sigma=F_s\left|_{s=0} \right. $.
Based on the results ( \ref{FE3}) and (\ref{MA11}), we deduce that $(\zeta,\sigma)$ is a solution to the problem
\begin{equation}
\begin{cases}
\begin{array}{c l}
\zeta_t+\sigma_{\xi}-\frac{\beta}{6}\zeta_{\xi\xi t}&=0 \\
\sigma_t+(c\zeta)_{\xi}-\frac{\beta}{6}\sigma_{\xi\xi t}&=((c-h)N)_{\xi},\\
\zeta(\xi,0)=\sigma(\xi,0)=0\\
\sigma(0,t)=\sigma(L,t)=0\\
\zeta(0,t)=\zeta(L,t)=0.
\end{array}
\end{cases}
\end{equation}
Furthermore, given that $c$ is solution to the RMP, we have
\begin{equation}
\frac{d}{ds}J_\alpha(c_s)\left|_{s=0}\right.\geq 0.
\end{equation}
Thus, from (\ref{MA13}), we deduce that
\begin{equation*}
\left\langle \zeta(\cdot,T) ,N(\cdot,T)-m_1(\cdot)\right\rangle_{H^{1}}\\+\left\langle \sigma(\cdot,T),V(\cdot,T)-m_2(\cdot) \right\rangle_{H^{1}}
+\alpha\left\langle c,(h-c)\right\rangle_{L^{2}}\geq 0,
\end{equation*}
for all $h \in \mathcal{M}_\gamma$. 

\noindent Next, multiplying the first equation of the adjoint problem (\ref{MA9}) by $\zeta$ and the second equation by $\sigma$, and considering $\eta_T(\xi)=N(\cdot,T;c)-m_1(\cdot)$, $\gamma_T(\xi)=V(\cdot,T;c)-m_2(\cdot)$ and integrating by parts over the interval $\left[0,L \right] \times \left[ 0,T\right] $, we obtain
\begin{equation}
\left\langle \zeta(\cdot,T),\eta_T(\cdot)\right\rangle_{H^{1}}+\left\langle \sigma(\cdot,T),\gamma_T(\cdot)\right\rangle_{H^{1}}=\left\langle\int_{0}^{T} N(\cdot,t)\gamma_\xi(\cdot,t)dt,(h-c)(\cdot)\right\rangle_{L^{2}}.
\end{equation}
Therefore, we conclude that
\begin{equation}
\left\langle\int_{0}^{T} N(\cdot,t)\gamma_\xi(\cdot,t)dt,(h-c)(\cdot)\right\rangle_{L^{2}}+\alpha\left\langle c,(h-c) \right\rangle_{L^{2}}\geq 0,
\end{equation}
which concludes the proof.
$\hfill\square$\\

\subsection{Analysis of an auxiliary problem}

In this section, we analyze an auxiliary problem associated with the adjoint problem in equation (\ref{MA9}).

Let $\tilde{c},c$ be coefficients in the space $L^{2}\left(\left[0,L\right] \right)$, and let $\tilde{U},U$ represent solutions to equation (\ref{De}) such that $\tilde{U}=(\tilde{N},\tilde{V})=U(\tilde{c})$ and $U=(N,V)=U(c)$. Define $\mathcal{N}=\tilde{N}-N$ and $\mathcal{V}=\tilde{V}-V$ as in equation (\ref{MA8}). Additionally, let $(\tilde{\eta},\tilde{\gamma})$ and $(\eta, \gamma)$ denote solutions to the adjoint problem (\ref{MA9}) associated with the coefficients $\tilde{c}$ and $c$, respectively, and assume they satisfy the final conditions
\begin{equation*}
\eta(\xi,T)=\eta_T(\xi), \\ \tilde{\eta}(\xi,T)=\tilde{\eta}_T(\xi)
\end{equation*}
\begin{equation*}
\gamma(\xi,T)=\gamma_T(\xi), \\ \tilde{\gamma}(\xi,T)=\tilde{\gamma}_T(\xi).
\end{equation*}
Define the functions $H(\xi,t;\tilde{c},c)$ and $R(\xi,t;\tilde{c},c)$ as
\begin{equation*}
H(\xi,t;\tilde{c},c)=H(\xi,t)=\tilde{\gamma}(\xi,t)-\gamma(\xi,t)
\end{equation*}
\begin{equation*}
R(\xi,t;\tilde{c},c)=R(\xi,t)=\tilde{\eta}(\xi,t)-\eta(\xi,t).
\end{equation*}
These functions, $H$ and $R$, belong to the space $\mathcal{L}^{2}_{T}$ and satisfy the problem
\begin{equation}\label{MA15}
\begin{cases}
\begin{array}{c l}
R_t+ cH_{\xi}-\frac{\beta}{6}R_{\xi\xi t}&=-(\tilde{c}-c)\tilde{\gamma}_\xi \\
H_t+R_{\xi}-\frac{\beta}{6}H_{\xi\xi t}&=0,\\
H(\xi,T)=H_T(\xi)=\tilde{\gamma}_T(\xi)-\gamma_T(\xi) \\
R(\xi,T)=R_T(\xi)=\tilde{\eta}_T(\xi)-\eta_T(\xi) \\
H(0,t)=H(L,t)=0\\
R(0,t)=R(L,t)=0.
\end{array}
\end{cases}
\end{equation}
The following result establishes an energy estimate for system (\ref{MA15}).

\begin{propo}
Given $R_T,H_T \in H^{1}_{0}(\left[0,L \right] )$, then the unique solution $(R,H)$ to problem (\ref{MA15}), satisfies the energy estimate
\begin{equation}\label{MA18}
\left\|(R,H) \right\|_{\mathcal{L}^{2}_{T}\times \mathcal{L}^{2}_{T}}\leq T^{\frac{1}{2}}\left[ \left\|(R_T,H_T) \right\|_{H^{1}\times H^{1}}+\frac{L^{\frac{1}{2}}}{(\frac{\beta}{6})}T^{\frac{1}{2}}\left\|\tilde{c}-c \right\|_{L^{2}}\left\|\tilde{\gamma} \right\|_{\mathcal{L}^{2}_{T}} \right] \text{exp} \left[ K_c T \right],
\end{equation}
where $K_c$ is a constant as in (\ref{MA7}).
\end{propo}
\noindent This result can be proved using the same strategy as in the proof of Proposition \ref{MA11}.

\subsection{Local stability and uniqueness of the RMP problem}\label{sec:Estabilidadyunicidadpi}

In this section, we establish a local stability estimate within the set of admissible coefficients for the restricted minimization problem (RMP), which is associated with the recovery of the linear velocity coefficient $c(x)$ in the linear Boussinesq system (\ref{De}) for initial data $(N_0,V_0)\in H_0^1 \times H_0^1$, based on final measurements $(m, \tilde{m}) \in H_0^1 \times H_0^1$, respectively, taken at the final time $T$. This constitutes the primary theoretical result of the present paper. It is worth noting that as a consequence, we deduce the uniqueness of the solution to the restricted minimization problem within this set of admissible coefficients.

\begin{teorema}\label{MA25}
Let $c,\tilde{c} \in \mathcal{M}_\gamma$ be optimal solutions to RMP associated with the measurements $m,\tilde{m}\in H^{1}_{0}(\left[0,L\right] )\times H^{1}_{0}(\left[0,L\right] )$, respectively. Then there exist $T_0> 0$ and a positive constant $\mathcal{S}$ such that for all $T\in \left. \left( 0,T_0\right]\right.$, we have 
\begin{equation}\label{MA19}
\left\|\tilde{c}-c\right\|_{L^{2}} \leq \mathcal{S} \left\|\tilde{m}-m \right\|_{H^{1}\times H^1} .
\end{equation}
\end{teorema}

\textit{Proof:} Let $m,\tilde{m}\in H^{1}_{0}(\left[0,L\right] )\times H^{1}_{0}(\left[0,L\right] )$, and let $c, \tilde{c} \in \mathcal{M}_\gamma$ be solutions to the RMP. Let $\tilde{U}=U(\tilde{c})=(\tilde{N},\tilde{V})$ and $U=U(c)=(N,V)$ be the solutions to the direct problem (\ref{De}), and let $(\tilde{\eta},\tilde{\gamma})$ and $(\eta,\gamma)$ be the corresponding solutions to the adjoint problem (\ref{MA9}) with final data
\begin{equation*}
\eta_T(\xi)=N(\xi,T)-m_1(\xi), \quad \quad \gamma_T(\xi)=V(\xi,T)-m_2(\xi),
\end{equation*}
and
\begin{equation*}
\tilde{\eta}_T(\xi)=\tilde{N}(\xi,T)-\tilde{m}_1(\xi), \quad \quad \tilde{\gamma}_T(\xi)=\tilde{V}(\xi,T)-\tilde{m}_2(\xi).
\end{equation*}
Observe that $c,\tilde{c}$ satisfy the optimality conditions (\ref{MA20}). Thus, for $h=\tilde{c}$, we have
\begin{equation*}
\left\langle \int^{T}_{0}N(\cdot,t)\gamma_\xi (\cdot,t)dt,(\tilde{c}-c)(\cdot)\right\rangle_{L^{2}}  +\alpha\left\langle c,(\tilde{c}-c)\right\rangle_{L^{2}}\geq 0.
\end{equation*}
Similarly, setting $h=c$, we obtain
\begin{equation*}
\left\langle \int^{T}_{0}\tilde{N}(\cdot,t)\tilde{\gamma}_\xi (\cdot,t)dt,(c-\tilde{c})(\cdot)\right\rangle_{L^{2}}  +\alpha\left\langle \tilde{c},(c-\tilde{c})\right\rangle_{L^{2}}\geq 0.
\end{equation*}
Adding these two inequalities gives
\begin{equation*}
\left\langle \int^{T}_{0}N\gamma_\xi-\tilde{N}\tilde{\gamma}_\xi dt,\tilde{c}-c\right\rangle_{L^{2}}  +\alpha\left\langle c-\tilde{c},\tilde{c}-c\right\rangle_{L^{2}}\geq 0.
\end{equation*}
Rearranging terms, we get
\begin{equation*}
\left\langle \int^{T}_{0}(\tilde{N}-N)\gamma_\xi+(\tilde{\gamma}_\xi-\gamma_\xi)\tilde{N} dt,c-\tilde{c}\right\rangle_{L^{2}}  -\alpha\left\|c-\tilde{c} \right\|^{2}_{L^{2}}\geq 0,
\end{equation*}
or equivalently,
\begin{equation}\label{MA21}
 \alpha\left\|c-\tilde{c} \right\|^{2}_{L^{2}}\leq \left\langle \int^{T}_{0}\mathcal{N}\gamma_\xi dt,c-\tilde{c}\right\rangle_{L^{2}} +\left\langle \int^{T}_{0}H_\xi\tilde{N} dt,c-\tilde{c}\right\rangle_{L^{2}},
\end{equation}
where $\mathcal{N}=\tilde{N}-N$ y $H=\tilde{\gamma}-\gamma$. 

\noindent Next, we estimate each term on the right-hand side of the previous inequality. We get
\begin{equation*}
\begin{split}
\left\langle \int^{T}_{0}\mathcal{N}\gamma_\xi dt,c-\tilde{c}\right\rangle_{L^{2}}&=\int_{0}^{T}\int_{0}^{L}(c-\tilde{c})\mathcal{N}\gamma_\xi d\xi dt\\
&\leq \int_{0}^{T}\left\|\mathcal{N}(t) \right\|_{L^{\infty}} \left\|c-\tilde{c} \right\|_{L^{2}}\left\|\gamma_\xi(t) \right\|_{L^{2}} dt\\
&\leq \frac{L^{\frac{1}{2}}}{\left( \frac{\beta}{6}\right) }\int_{0}^{T}\left\|\mathcal{N}(t) \right\|_{H^{1}} \left\|c-\tilde{c} \right\|_{L^{2}}\left\|\gamma(t) \right\|_{H^{1}} dt\\
&\leq \frac{L^{\frac{1}{2}}}{2\left( \frac{\beta}{6}\right) }\int_{0}^{T}\left( \frac{1}{T}\left\|\mathcal{N}(t) \right\|_{H^{1}}^{2} +T\left\|c-\tilde{c} \right\|_{L^{2}}^{2}\left\|\gamma(t) \right\|_{H^{1}}^{2} \right) dt\\
&\leq \frac{L^{\frac{1}{2}}}{2\left( \frac{\beta}{6}\right) }\left( \frac{1}{T}\left\|\mathcal{N}\right\|_{\mathcal{L}^{2}_{T}}^{2} +T\left\|c-\tilde{c} \right\|_{L^{2}}^{2}\left\|\gamma \right\|_{\mathcal{L}^{2}_{T}}^{2} \right). \\
\end{split}
\end{equation*}
Using the energy estimates (\ref{MA8}) and (\ref{MA9}), we arrive at
\begin{multline}\label{MA23}
\left\langle \int^{T}_{0}\mathcal{N}\gamma_\xi dt,c-\tilde{c}\right\rangle_{L^{2}}
\leq \frac{L^{\frac{3}{2}}T^{2}}{2\left( \frac{\beta}{6}\right)^{3} }\left\|c-\tilde{c} \right\|_{L^{2}}^{2}\left\|(N_0,V_0) \right\|_{H^{1}\times H^{1}}^{2}K_1(c,\tilde{c})^{2}\\+\frac{L^{\frac{1}{2}}T^{2}}{2\left( \frac{\beta}{6}\right) }\left\|c-\tilde{c} \right\|_{L^{2}}^{2}\left\|(\eta_T,\gamma_T) \right\|_{H^{1}\times H^{1}}^{2}K_2(c)^{2},
\end{multline}
where $K_2(c)$ y $K_1(c,\tilde{c})$ are as in (\ref{FE1}), (\ref{MA8}), respectively. \\

\noindent On the other hand, the last term in (\ref{MA21}) can be estimated as
\begin{equation*}
\begin{split}
\left\langle \int^{T}_{0}H_\xi\tilde{N} dt,c-\tilde{c}\right\rangle_{L^{2}}&=\int_{0}^{T}\int_{0}^{L}(c-\tilde{c})H_\xi\tilde{N} d\xi dt\\
&\leq \int_{0}^{T}\left\|\tilde{N}(t) \right\|_{L^{\infty}} \left\|c-\tilde{c} \right\|_{L^{2}}\left\|H_\xi(t) \right\|_{L^{2}} dt\\
&\leq \frac{L^{\frac{1}{2}}}{\left( \frac{\beta}{6}\right) }\int_{0}^{T}\left\|\tilde{N}(t) \right\|_{H^{1}} \left\|c-\tilde{c} \right\|_{L^{2}}\left\|H(t) \right\|_{H^{1}} dt\\
&\leq \frac{L^{\frac{1}{2}}}{2\left( \frac{\beta}{6}\right) }\int_{0}^{T}\left( \frac{1}{T}\left\|H(t) \right\|_{H^{1}}^{2} +T\left\|c-\tilde{c} \right\|_{L^{2}}^{2}\left\|\tilde{N}(t) \right\|_{H^{1}}^{2} \right) dt\\
&\leq \frac{L^{\frac{1}{2}}}{2\left( \frac{\beta}{6}\right) }\left( \frac{1}{T}\left\|H \right\|_{\mathcal{L}^{2}_{T}}^{2} +T\left\|c-\tilde{c} \right\|_{L^{2}}^{2}\left\|\tilde{N} \right\|_{\mathcal{L}^{2}_{T}}^{2} \right). \\
\end{split}
\end{equation*}
Consequently, using the energy estimates (\ref{MA18}) and (\ref{Ana}), we obtain
\begin{multline}
\left\langle \int^{T}_{0}H_\xi \tilde{N} dt,c-\tilde{c}\right\rangle_{L^{2}}
\leq \frac{L^{\frac{1}{2}}T^{2}}{2\left( \frac{\beta}{6}\right) }\left\|c-\tilde{c} \right\|_{L^{2}}^{2}\left\|(N_0,V_0) \right\|_{H^{1}\times H^{1}}^{2}K_2(\tilde{c})^{2}\\+\frac{L^{\frac{1}{2}}}{\left( \frac{\beta}{6}\right) }\left\|(R_T,H_T) \right\|_{H^{1}\times H^{1}}^{2}K_2(c)^{2}+\frac{2L^{\frac{3}{2}}T}{\left( \frac{\beta}{6}\right)^{3} }\left\|c-\tilde{c} \right\|_{L^{2}}^{2}\left\|\tilde{\gamma} \right\|^{2}_{\mathcal{L}^{2}_{T}}K_2(c)^{2}.
\end{multline}
Note that
\begin{equation*}
\left\|H_T\right\|^{2}_{H^{1}}\leq 2\left( \left\| \left( \tilde{V}-V\right)\left(\cdot, T\right)  \right\|^{2}_{H^{1}} +\left\|\tilde{m}_2 -m_2\right\|^{2}_{H^{1}} \right),
\end{equation*}
\text{and}
\begin{equation*}
\left\|R_T\right\|^{2}_{H^{1}}\leq 2\left( \left\| \left( \tilde{N}-N\right)\left(\cdot, T\right)  \right\|^{2}_{H^{1}} +\left\|\tilde{m}_1 -m_1\right\|^{2}_{H^{1}} \right).
\end{equation*}
Therefore, applying the energy estimates (\ref{MA8}), (\ref{MA9}), (\ref{MA18}), we get
\begin{multline}\label{MA22}
\left\langle \int^{T}_{0}H_\xi \tilde{N} dt,c-\tilde{c}\right\rangle_{L^{2}}
\leq  \left[ \frac{L^{\frac{1}{2}}T^{2}}{2\left( \frac{\beta}{6}\right) }\left\|(N_0,V_0) \right\|_{H^{1}\times H^{1}}^{2}K_2(\tilde{c})^{2}\right. \\ \left. \frac{2L^{\frac{3}{2}}T^{2}}{\left( \frac{\beta}{6}\right)^{3} }\left\|(\tilde{\eta}_T,\tilde{\gamma}_T) \right\|^{2}_{H^{1} \times H^{1}}K_2(\tilde{c})^{2}K_2(c)^{2}+ \frac{8L^{\frac{3}{2}}T^{3}}{\left( \frac{\beta}{6}\right)^{3} }\left\|(N_0,V_0) \right\|_{H^{1}\times H^{1}}^{2}K_1(c,\tilde{c})^{2}K_2(c)^{2}\right]\left\|c-\tilde{c} \right\|^{2}_{L^{2}}\\
+\frac{4L^{\frac{1}{2}}}{\left( \frac{\beta}{6}\right) }\left[\left\|\tilde{m}_2 -m_2\right\|^{2}_{H^{1}}+\left\|\tilde{m}_1 -m_1\right\|^{2}_{H^{1}} \right]K_2(c)^{2}.
\end{multline}
Now, substituting (\ref{MA23}), (\ref{MA22}) into (\ref{MA21}), we obtain
\begin{equation}
\alpha \left\|c-\tilde{c} \right\|^{2}_{L^{2}}\leq \mathcal{H}(T)\left\|c-\tilde{c} \right\|^{2}_{L^{2}}+\mathcal{W}\left[\left\|\tilde{m}_2 -m_2\right\|^{2}_{H^{1}}+\left\|\tilde{m}_1 -m_1\right\|^{2}_{H^{1}}\right]
\end{equation}
\begin{equation}
\alpha \left\|c-\tilde{c} \right\|^{2}_{L^{2}}\leq \mathcal{H}(T)\left\|c-\tilde{c} \right\|^{2}_{L^{2}}+\mathcal{W}\left\|\tilde{m} -m\right\|^{2}_{H^{1}},
\end{equation}
where
\begin{multline}\label{MA24}
\mathcal{H}(T)= \frac{L^{\frac{1}{2}}T^{2}}{2\left( \frac{\beta}{6}\right) }\left\|(N_0,V_0) \right\|_{H^{1}\times H^{1}}^{2}K_2(\tilde{c})^{2}  \frac{2L^{\frac{3}{2}}T^{2}}{\left( \frac{\beta}{6}\right)^{3} }\left\|(\tilde{\eta}_T,\tilde{\gamma}_T) \right\|^{2}_{H^{1} \times H^{1}}K_2(\tilde{c})^{2}K_2(c)^{2}\\+ \frac{8L^{\frac{3}{2}}T^{3}}{\left( \frac{\beta}{6}\right)^{3} }\left\|(N_0,V_0) \right\|_{H^{1}\times H^{1}}^{2}K_1(c,\tilde{c})^{2}K_2(c)^{2}+\frac{L^{\frac{3}{2}}T^{2}}{2\left( \frac{\beta}{6}\right)^{3} }\left\|(N_0,V_0) \right\|_{H^{1}\times H^{1}}^{2}K_1(c,\tilde{c})^{2}\\+\frac{L^{\frac{1}{2}}T^{2}}{2\left( \frac{\beta}{6}\right) }\left\|(\eta_T,\gamma_T) \right\|_{H^{1}\times H^{1}}^{2}K_2(c)^{2},
\end{multline}
and
\begin{equation}
\mathcal{W}=\frac{4L^{\frac{1}{2}}}{\left( \frac{\beta}{6}\right) }K_2(c)^{2}.
\end{equation}

\noindent Finally, selecting $T_0>0$ in such way that $0<\mathcal{H}(T_0)<\alpha$, we can conclude that
\begin{equation}
\left\|c-\tilde{c} \right\|^{2}_{L^{2}}\leq \frac{\mathcal{W}}{\alpha-\mathcal{H}(T)}\left\|\tilde{m} -m\right\|^{2}_{H^{1} \times H^1},
\end{equation}
for any $0<T<T_0$, which concludes the proof.

$\hfill\square$


Theorem (\ref{MA25}) implies the following result.\\

\begin{corolario}
Suppose the hypotheses in theorem (\ref{MA25}). If
\begin{equation*}
m(\xi)=\tilde{m}(\xi), \quad \textup{para todo} \quad \xi \in \left[ 0,L\right],
\end{equation*}
then there exists $T_0>0$, such as for all $T \in \left(0,T_0\right] $, 
\begin{equation*}
c(\xi)=\tilde{c}(\xi),\quad \textup{ for almost all }\quad \xi \in \left[ 0,L\right]. 
\end{equation*}
\end{corolario}

\section{Numerical Simulations}	

Let $0< L_0 < L_1$. It is important to note that the theoretical results presented in the previous section were derived under the assumption of linearization given by \eqref{KdVe1}. In this section, we employ a numerical approach to investigate the retrieval of the first-order coefficient $M(\xi)$ in the fully nonlinear-dispersive system \eqref{KdVe}, which extends beyond the theoretical findings discussed in this paper.

Firstly, we introduce a suitable numerical scheme to approximate the solution of the direct problem associated with the Boussinesq system
\begin{equation}\label{MA29}
  \begin{array}{c l}
   M(\xi) \eta_t+\partial _\xi \left[\left(1+\frac{ \tilde{\alpha} \eta}{M\left(\xi\right)}\right)u\right]-\frac{\beta}{6}\partial^2 _\xi \left(M(\xi) \eta_t\right)&=0, \:\: \left(\xi,t\right)\in [L_0,L_1] \times \left[0, T \right], \\
  u_t+\eta_\xi+\frac{ \tilde{\alpha} }{2}\partial _\xi\left[\left(\frac{u}{M\left(\xi\right)}\right)^2\right]-\frac{\beta}{6}\partial^2 _\xi \left(u_t\right)&=0.
\end{array}
\end{equation}
subject to the initial conditions 
\begin{equation}\label{datosiniciales2}
\begin{array}{c l}
\eta(\xi,0)&=\eta_0(\xi),\\
u(\xi,0)&=u_0(\xi),\:\:\: \xi \in \left[L_0,L_1\right],
\end{array}
\end{equation}
and Dirichlet boundary conditions
\begin{equation}\label{condfrontera2}
\begin{array}{c l}
\eta(L_0,t)=\eta(L_1,t)=0,\:\: t\in \left[0,T\right]\\
u(L_0,t)=u(L_1,t)=0,\:\: t\in \left[0,T\right].
\end{array}
\end{equation}

\subsection{Numerical scheme for the full nonlinear Boussinesq system}

In this section, we propose a finite element numerical scheme in the spatial variable to approximate solutions of the nonlinear Boussinesq system 
(\ref{MA29}). To construct this scheme, we first derive the variational formulation of this problem. Specifically, we multiply the first equation of system \eqref{MA29} by
a test function $v_1$, and the second equation by another test function $v_2$, both chosen from the Sobolev space $H^{1}_{0}$.  Integrating the resulting equations over the interval $\left[L_0,L_1 \right] $, and applying integration by parts, we obtain the system

\begin{equation*}
\left\langle M \eta_t , v_1 \right\rangle_{L^2} - \left\langle \Big(1 + \frac{\tilde{\alpha} \eta}{M} \Big) u , \partial_\xi v_1 \right\rangle_{L^2} + \frac{\beta}{6} \left\langle \partial_\xi(M \eta_t), \partial_\xi v_1 \right\rangle_{L^2}  =  0,
\end{equation*}

\begin{equation*}
\left\langle u_t, v_2 \right\rangle_{L^2} - \left\langle \eta, \partial_\xi v_2 \right\rangle_{L^2} - \frac{\tilde{\alpha}}{2} \left\langle \Big( \frac{u}{M} \Big)^2 , \partial_\xi v_2 \right\rangle_{L^2} +
\frac{\beta}{6} \left\langle \partial_\xi u_t, \partial_\xi v_2 \right\rangle_{L^2}  =  0.
\end{equation*}
Observe that the previous equations can be cast into the system

\begin{equation}\label{sistema1}
\partial_t\left(    \left\langle M \eta,v_1\right\rangle _{L^{2}} + \frac{\beta}{6}\left\langle \partial_\xi (M \eta ),\partial_\xi v_1 \right\rangle _{L^{2}} \right)   =  \left\langle \Big(1 + \frac{\tilde{\alpha} \eta}{M} \Big) u , \partial_\xi v_1 \right\rangle_{L^2}, 
\end{equation}

\begin{equation}\label{sistema2}
\partial_t\left( \left\langle u,v_2\right\rangle _{L^{2}}+\frac{\beta}{6}\left\langle \partial_\xi u, \partial_\xi v_2 \right\rangle _{L^{2}}\right) =    \left\langle \eta, \partial_\xi v_2 \right\rangle_{L^2} + \frac{\tilde{\alpha}}{2} \left\langle \Big( \frac{u}{M} \Big)^2 , \partial_\xi v_2 \right\rangle_{L^2},
\end{equation}
which in turn can be viewed as a system of ordinary differential equations in the form
\begin{equation}\label{EDOBoussinesq}
U_t=F(U),
\end{equation}
where $U = (U_1, U_2)$, with
\begin{align*}
&U_1 := \left\langle M \eta,v_1\right\rangle _{L^{2}} + \frac{\beta}{6}\left\langle \partial_\xi (M \eta ),\partial_\xi v_1 \right\rangle _{L^{2}}  \\
&U_2 := \left\langle u,v_2\right\rangle _{L^{2}}+\frac{\beta}{6}\left\langle \partial_\xi u, \partial_\xi v_2 \right\rangle _{L^{2}},
\end{align*}
and $F(U) := (F_1(U), F_2(U) )$, with
\begin{align*}
&F_1(U ) :=  \left\langle \Big(1 + \frac{\tilde{\alpha} \eta}{M} \Big) u , \partial_\xi v_1 \right\rangle_{L^2},\\
&F_2(U) :=  \left\langle \eta, \partial_\xi v_2 \right\rangle_{L^2} + \frac{\tilde{\alpha}}{2} \left\langle \Big( \frac{u}{M} \Big)^2 , \partial_\xi v_2 \right\rangle_{L^2}.
\end{align*}
Equation \eqref{EDOBoussinesq} is discretized using an implicit second-order scheme, expressed as
\begin{equation}\label{ApproxBoussinesq}
\frac{U^{n+1}-U^{n}}{\Delta t}=\theta F(U^{n+1})+(1-\theta)F(U^{n}), 
\end{equation}
where $\theta = 1/2$, and $U^n$ represents the approximation of the function $U$ at the time $t=n\Delta t$.

The integrals with respect to the spatial variable $\xi$ in the formulation \eqref{EDOBoussinesq} are discretized using the Finite Element Method (FEM) on the interval $[L_0,L_1]$, as implemented by the Python libraries Python-FEniCS-Dolfin (\cite{Fenics, Logg, Fenics2, Farrell, Funke, Alnaes1, Alnaes1b, Alnaes2}). To apply the FEM, the interval $[L_0,L_1]$ is divided into a regular partition $(\xi_j)_{1 \leq j \leq N}$, with step size$\Delta \xi = \xi_j - \xi_{j-1}$, 
and the solution $(\eta(\xi,t),u(\xi,t))$ at time $t = n \Delta t$ is approximated by
\begin{equation*}
\eta^{n}(\xi)=\sum_{i=1}^{N}\eta_i^{n}\phi_i(\xi)
\end{equation*}
and
\begin{equation*}
u^{n}(\xi)=\sum_{i=1}^{N}u_i^{n}\phi_i(\xi),
\end{equation*}
respectively, where $\phi_i(\xi)$, $i=1,...,N$ denotes a set of piecewise linear functions defined on the interval $[L_0,L_1]$,
and $\eta_i^n, u_i^n$ are constant coefficients.

The approach described above, where the temporal variable is discretized first, followed by the spatial variable, is known as Rothe's method \cite{Kacur},\cite{Rothe}. This strategy is similar to the methodology employed by the authors in  \cite{JCPipicanoSosa} for a linear BBM equation, and we have adapted it to the general nonlinear-dispersive Boussinesq system \eqref{MA29}. 

\subsection{Approximation of the restricted minimization problem}
To numerically explore the inverse problem of identifying the coefficient $M(\xi)$ in system \eqref{MA29} over the interval domain $[L_0, L_1]$, we formulate a restricted minimization problem (RMP) as follows:
\[
\min_{M \in \mathcal{A}_\epsilon } J_{\alpha, 2} ( M ),
\]
where $J_{\alpha, 2}$ is defined by
\begin{align}\label{functional2a}
&J_{\alpha, 2}(M) := \frac12 \int_{L_0}^{L_1} | \eta(\xi, T, M) - \eta_T(\xi) |^2 d\xi +\frac12 \int_{L_0}^{L_1} | u(\xi, T, M) - u_T(\xi) |^2 d\xi + \\
&~~~~~~~~~~~~ \frac{\alpha}{2} \int_{L_0}^{L_1} | M(\xi) -1 |^2 d\xi,  \notag
\end{align}
which serves as a regularized version of the original functional $J_{\alpha}$. The original $J_\alpha$ led to numerical oscillations in the computed coefficients due to the derivative terms involving the dispersion parameter $\beta$.
The regularization term used here normalizes the coefficient around the constant reference value 1. As demonstrated in the numerical experiments, this term serves to reinforce the coefficient's alignment with the
reference value within the numerical solution space. 

Given observations data $\eta_T, u_T$, and a discrete coefficient $M =(M_j)_{1 \leq j \leq N}$ measured at the mesh points $(\xi_j)_{j=1,...,N}$ along the interval $[L_0,L_1]$, we consider the discrete version of the functional $J_{\alpha,2}$ as
\begin{equation}
\begin{aligned}\label{functional2b}
&\mathcal{J}_{\alpha,2} (M) := \frac12 \sum_{j=1}^N  | \eta(\xi_j, T, M) - \eta_T(\xi_j) |^2 \Delta \xi   + \frac12 \sum_{j=1}^N  | u(\xi_j, T, M) - u_T( \xi_j) |^2 \Delta \xi  + \\
& ~~~~~~~~~~~~~\frac{\alpha}{2} \sum_{j=1}^N |M_j-1|^2 \Delta \xi, 
\end{aligned}
\end{equation}
where $(\eta(\xi, T, M), u(\xi, T, M) )$ represents the approximation of the solution to problem \eqref{MA29} at time $t = T$ for the discrete coefficient $M$. This functional quantifies the discrepancy between the observed data and the solution of problem \eqref{MA29} at the final time $T$. Minimization of $\mathcal{J}_\alpha$ is performed using the Dolfin-Adjoint library \cite{Farrell, Funke}, via the iterative L-BFGS-B algorithm. This algorithm has been successfully employed in previous studies by one of the authors in \cite{JCPipicanoSosa} for an inverse problem involving the linear BBM equation. This algorithm is implemented in the SciPy library \cite{scipy}, as described by \cite{Byrd, Zhu}.

To ensure convergence, the parameter \text{ftol} in the L-BFG-S routine halts the iteration when consecutive values of the objective functional satisfy the condition

\begin{equation}\label{stop_condition}
\frac{| \mathcal{J}_{\alpha,2}^{(k+1)} - \mathcal{J}_{\alpha,2}^{(k)} |}{ \max \{1, | \mathcal{J}_{\alpha,2}^{(k+1)}|, |\mathcal{J}_{\alpha,2}^{(k)}| \} } \leq \text{ftol}.
\end{equation}

In the numerical simulations presented, unless otherwise stated, we set $\text{ftol} = 10^{-8}$. Here $\mathcal{J}_\alpha^{(k)}$ denotes the value of the discrete objective functional $\mathcal{J}_\alpha$ evaluated at the coefficient $M^{(k)}$, which is computed at the $k$-th iteration of the L-BFGS-B algorithm.

\subsection{Description of the numerical experiments}

Unless explicitly stated otherwise, we consider the initial condition for the nonlinear system \eqref{MA29} to be the smooth compactly supported profile:
$$
\eta(\xi,0) = \eta_0(\xi) = u(\xi,0) =  u_0(\xi) = e^{-(\xi - 3)^2},
$$
where the computational domain is the interval $[L_0, L_1]=[-20,40]$. Additionally, the initial guess for the unknown coefficient $M$ in the minimization process of the functional $\mathcal{J}_{\alpha,2}$ is set to $M \equiv 1$, representing a homogeneous medium. All absolute errors are computed using the $L^2$-norm.

To demonstrate the efficiency of the proposed numerical strategy for approximating the solution of the inverse problem,
where we seek to reconstruct the coefficient $M(\xi)$ in the nonlinear system \eqref{MA29} by solving the restricted minimization problem (RMP), we approximate the solutions $\eta(\xi,t, M)$ and $u(\xi,t, M)$ of the direct problem \eqref{MA29}, using the numerical scheme provided in \eqref{ApproxBoussinesq}. The resulting profiles $\eta(\xi, T, M)$ and $u(\xi, T, M)$ at time $t = T$ serve as our synthetic target functions for the reconstruction process of the coefficient $M(\xi)$, employing the L-BFGS-B algorithm mentioned earlier with a stopping parameter $\text{ftol} = 10^{-8}$.

\vspace{0.5cm}

{\bf Experiment set 1:} Recovery of a smooth coefficient.

\begin{figure}
\centering
\includegraphics[width=10cm, height=8cm]{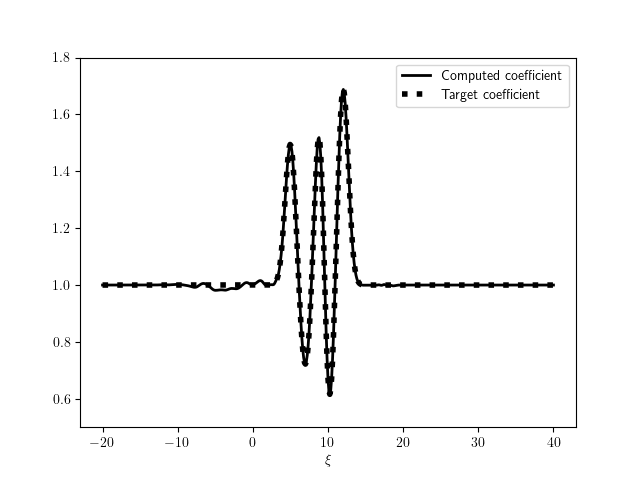}
\caption{Computed coefficient $M(\xi)$ after 25 iterations in Experiment set 1, with parameters $\tilde{\alpha}=0$, $\beta=0.1$, and $\alpha=0$ (no regularization term). }
\label{Exp1NEW}
\end{figure}

In this initial experiment, we investigate the nonlinear system \eqref{MA29} with the following parameters: $\beta = 0.1$, $\tilde{\alpha} = 0$ (i.e., linear regime), $T = 15$, $\Delta x = 60/500 = 0.12$, and $\Delta t = T/1500 \approx 0.01$. Additionally, we define the exact coefficient $M(\xi)$ as the smooth function

\begin{equation}\label{gauss_coeff}
M(\xi) = 1 + 0.5 e^{-(\xi - 5)^2} -0.3 e^{ - (\xi - 7)^2} + 0.7 e^{-(\xi - 9)^2} - 0.6 e^{-(\xi - 10 )^2} + 
0.7 e^{-(\xi - 12)^2},
\end{equation}
and $\alpha = 0$, indicating no regularization in the objective functional
$\mathcal{J}_{\alpha,2}$. The result of the simulation is presented in Figure \ref{Exp1NEW}. After 25 iterations of the L-BFGS-B minimization method, the absolute error between the exact and computed coefficients is approximately 0.036.

\begin{figure}
\centering
\includegraphics[width=10cm, height=8cm]{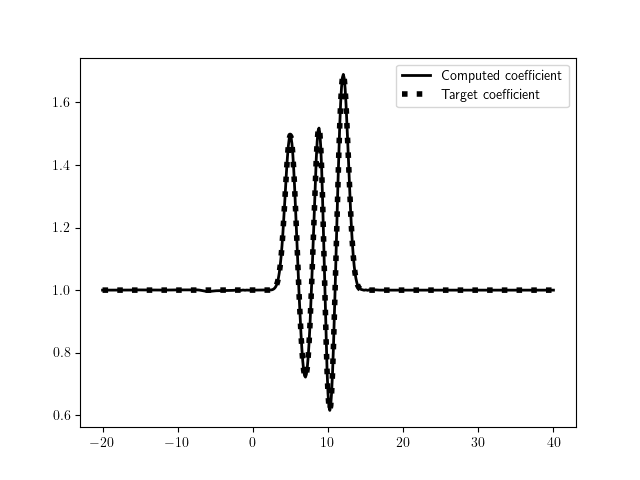}
\caption{Computed coefficient $M(\xi)$ after 130 iterations in Experiment set 1, with parameters $\tilde{\alpha}=0$, $\beta=0.1$, and $\alpha=0$ (no regularization term). }
\label{Exp1_130it}
\end{figure}

In contrast, the result with 130 iterations is shown in Figure \ref{Exp1_130it}, where an absolute error of approximately $5 \times 10^{-3}$ is observed between the numerical and exact coefficients. We observe that the number of iterations has a hidden regularization effect, since in these experiments we have not used directly a penalty term in the functional $\mathcal{J}_{\alpha,2}$.

\vspace{0.5cm}

{\bf Experiment set 2:} Recovery of a non-smooth coefficient.

In second experiment, we maintain the same model and numerical parameters as in the previous one, but with a different
final time is $T = 20$ (resulting in $\Delta t \approx T/1500 = 0.133$), and initial data given by
$$
\eta_0(\xi) = u_0(\xi) = e^{- \xi^2}.
$$
The exact coefficient is now represented by the piecewise constant function
\begin{equation}\label{irregularcoeff1}
M(\xi) =  \begin{cases}
               &1.3, ~~  5 \leq \xi < 15\\
               &1, ~~ \text{otherwise}.
           \end{cases}
\end{equation}

\begin{figure}
\centering
\includegraphics[width=10cm, height=8cm]{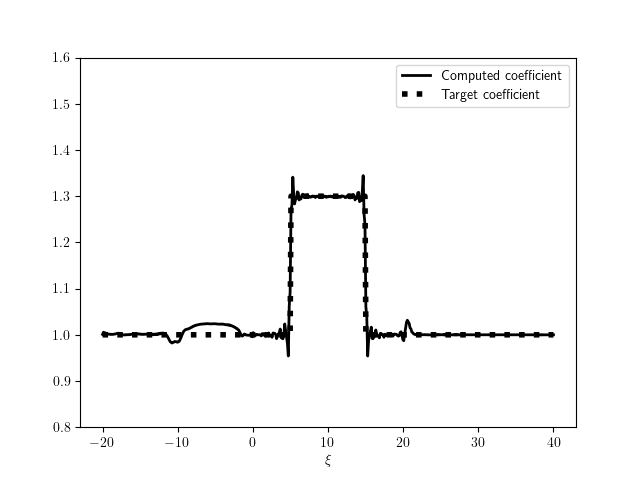}
\caption{Computed coefficient $M(\xi)$ after 62 iterations in Experiment set 2, with parameters $\tilde{\alpha}=0$, $\beta=0.1$, and $\alpha=0$ (no regularization term). }
\label{Exp2NEW}
\end{figure}

\begin{figure}
\centering
\includegraphics[width=10cm, height=8cm]{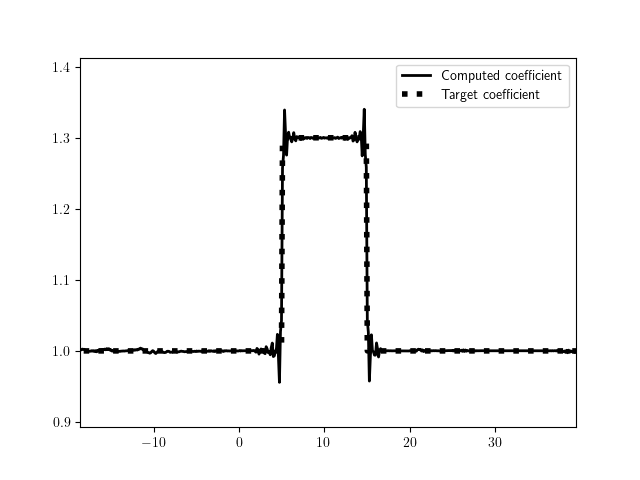}
\caption{Computed coefficient $M(\xi)$ after 500 iterations in Experiment set 2, with parameters $\tilde{\alpha}=0$, $\beta=0.1$, and regularization parameter $\alpha=10^{-5}$. }
\label{Coeff2_with_regul_cuad}
\end{figure}

We set $\alpha=0$, indicating that no regularization in the objective functional $\mathcal{J}_{\alpha, 2}$. The result of this computer simulation is shown in Figure \ref{Exp2NEW}. After 62 iterations of the L-BFGS-B minimization method, the absolute $L^2$-error between the exact and computed coefficients is approximately 0.0826. In contrast, Figure \ref{Coeff2_with_regul_cuad} presents the results of a similar simulation using the same numerical and modeling parameters, but now incorporating a regularization term with $\alpha = 10^{-5}$ in the functional $\mathcal{J}_{\alpha, 2}$. After 500 iterations, the absolute $L^2$-error between computed and exact coefficient is reduced to approximately 0.04 - half the error obtained without regularization.

These results demonstrate that incorporating regularization significantly improves the accuracy of the coefficient reconstruction process, especially in the case of a non-differentiable coefficient.

\vspace{0.5cm}






{\bf Experiment set 3:} Effect of the nonlinear parameter $\tilde{\alpha}$ in coefficient identification.

In this set of numerical experiments, we investigate the performance of the coefficient identification scheme in the presence of nonlinear and dispersive terms in the Boussinesq system \eqref{MA29}. We fix the dispersive parameter at $\beta = 0.1$ and progressively increase the value of the nonlinear parameter $\tilde{\alpha}$ to assess its influence on the coefficient reconstruction process.
The target coefficient $M(\xi)$ is modeled as the piecewise constant function \eqref{irregularcoeff1},
previously considered in Experiment set 2.

\begin{figure}[ht]
  \centering
	   \begin{subfigure}{0.35\linewidth}
		\includegraphics[width=\linewidth, height=\linewidth]{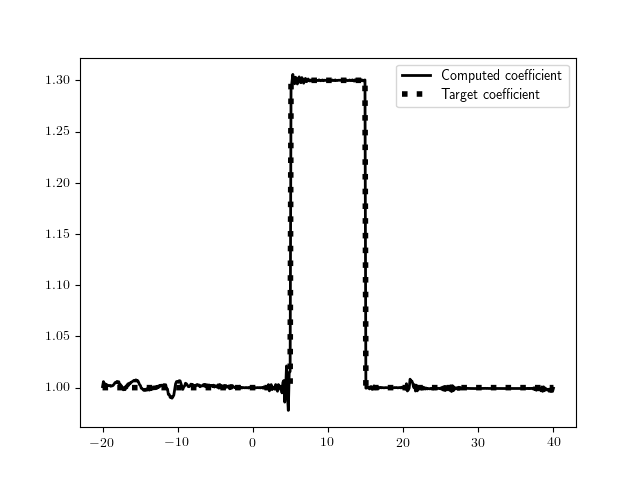}
		\caption{$\tilde{\alpha} = 0.01$. }
		\label{eff_nonlin1}
	   \end{subfigure}
	   \begin{subfigure}{0.35\linewidth}
	      	\includegraphics[width=\linewidth, height=\linewidth]{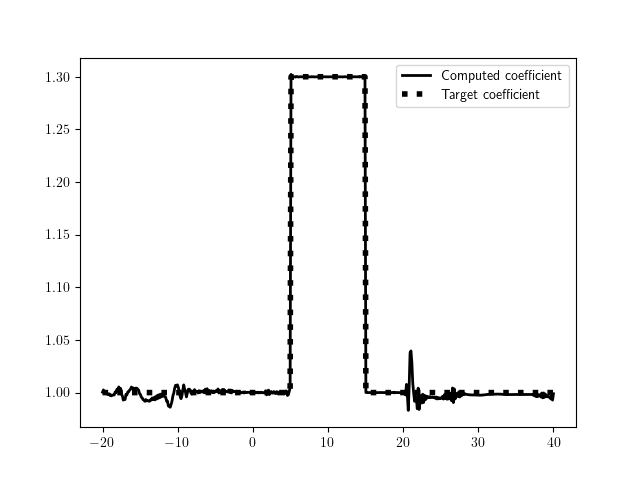}
		\caption{$\tilde{\alpha} = 0.03$. }
		\label{eff_nonlin2 }
	    \end{subfigure}
	\vfill
	     \begin{subfigure}{0.35\linewidth}
		 \includegraphics[width=\linewidth, height=\linewidth]{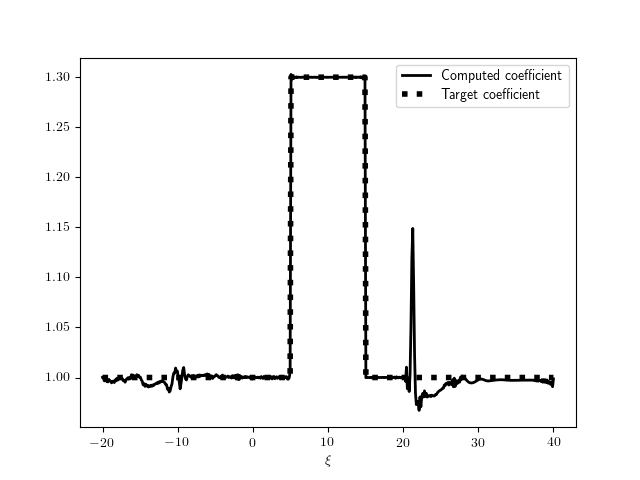}
		 \caption{$\tilde{\alpha} = 0.05$. }
		 \label{eff_nonlin3}
	      \end{subfigure} 
	       \begin{subfigure}{0.35\linewidth}
		  \includegraphics[width=\linewidth, height=\linewidth ]{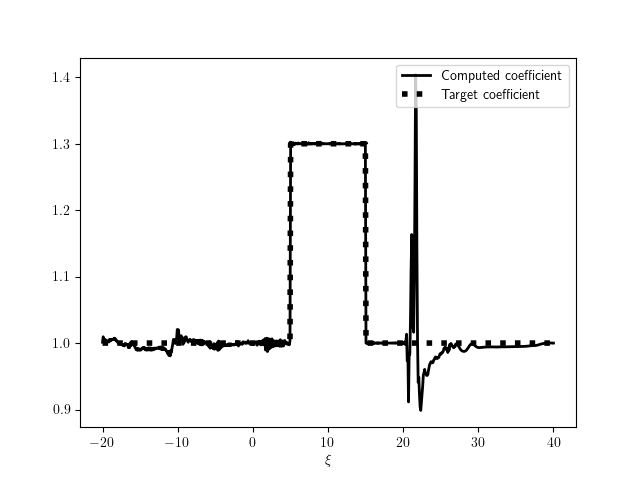}
		  \caption{$\tilde{\alpha}=0.07$.}
		  \label{eff_nonlin4}
	       \end{subfigure}   
	\caption{Computation of the coefficient $M(\xi)$ in Experiment set 3, with $\beta=0.1$, for various values of the nonlinear parameter $\tilde{\alpha}$, and without a regularizing term (i.e. $\alpha = 0)$. The absolute $L^2$-errors are as follows : (a) 0.02, (b) 0.032, (c) 0.086,  (d) 0.21. }
	\label{effect_nonlinearity}
\end{figure}

The results of the simulations are presented in Figure \ref{effect_nonlinearity}, after 500 iterations of the L-BFGS-B method. For subfigures (a), (b), and (c), the numerical parameters are identical to those in Experiment set 2. For subfigure (d), however, the parameters were adjusted to $\Delta t = 20/1700 \approx 0.012$ and $\Delta x = 60/700 \approx 0.086$ to ensure convergence of the Newton solver employed for the Boussinesq system \eqref{MA29}. 

As the nonlinear parameter $\tilde{\alpha}$ increases, the optimization algorithm requires more iterations to approximate the target coefficient, and the resulting $L^2$-error increases accordingly.
Notably,  for $\tilde{\alpha}=0.05$ and $\tilde{\alpha}=0.07$, the $L^2$-errors reach 0.086 and 0.21, respectively.

\begin{figure}[ht]
  \centering
	   \begin{subfigure}{0.35\linewidth}
		\includegraphics[width=\linewidth, height=\linewidth]{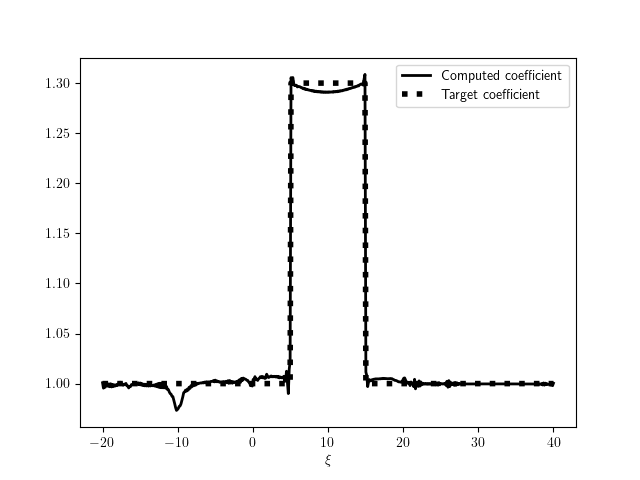}
		\caption{$\tilde{\alpha} = 0.05$. }
		\label{eff_nonlinnew1}
	   \end{subfigure}
	   \begin{subfigure}{0.35\linewidth}
	      	\includegraphics[width=\linewidth, height=\linewidth]{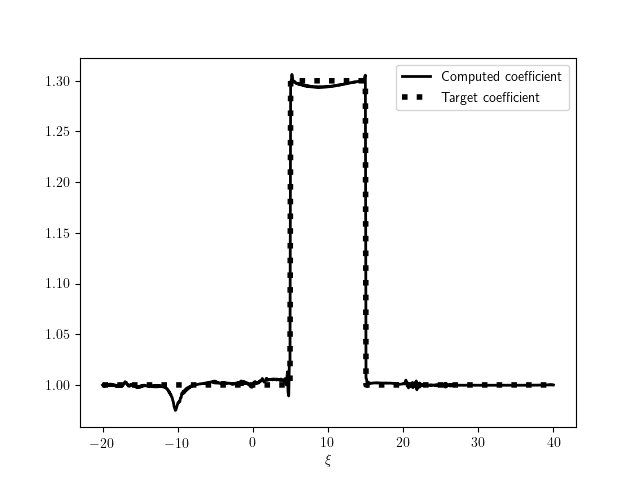}
		\caption{$\tilde{\alpha} = 0.07$. }
		\label{eff_nonlinnew2 }
	    \end{subfigure}
	\caption{Computation of the coefficient $M(\xi)$ in Experiment set 3, with $\beta=0.1$, for $\tilde{\alpha}=0.05$ (a) and $\tilde{\alpha}=0.07$ (b), using the functional $\mathcal{J}_{\alpha,2}$ with $\alpha = 0.001$. The absolute $L^2$-errors result to be approximately 0.04 and 0.03, respectively. }
	\label{effect_nonlinearitynew2}
\end{figure}

To mitigate this degradation in accuracy, we employ the regularized functional $\mathcal{J}_{\alpha,2}$ with a regularization parameter $\alpha = 0.001$. Using this approach, the results after 500 iterations show significant improvement. In Figure \ref{effect_nonlinearitynew2}(a), with $\tilde{\alpha}=0.05$, $\beta=0.1$, $\Delta t = 20/1500$, and $\Delta x = 60/500$, the $L^2$-error is reduced to approximately 0.04. Similarly, in Figure \ref{effect_nonlinearitynew2}(b), for $\tilde{\alpha}=0.07$, $\Delta t = 20/1700 \approx 0.0117$, and $\Delta x = 60/700 \approx 0.0857$, the error decreases to about 0.03.

These experiments clearly demonstrate that the inclusion of a regularization term in the Tikhonov functional becomes increasingly important as the nonlinear parameter $\tilde{\alpha}$ grows.



\begin{figure}[ht]
  \centering
	   \begin{subfigure}{0.35\linewidth}
		\includegraphics[width=\linewidth, height=\linewidth]{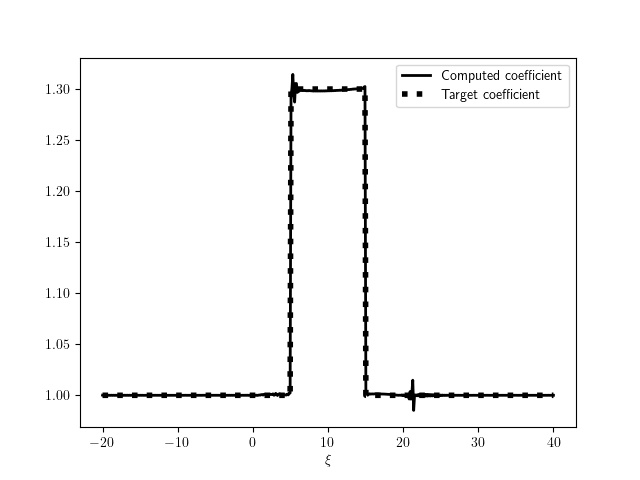}
		\caption{$\tilde{\alpha} = 0.05$. }
		\label{eff_nonlin1}
	   \end{subfigure}
	   \begin{subfigure}{0.35\linewidth}
	      	\includegraphics[width=\linewidth, height=\linewidth]{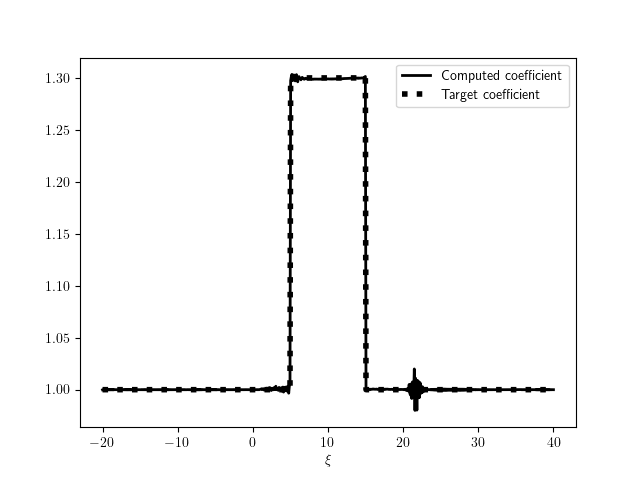}
		\caption{$\tilde{\alpha} = 0.07$. }
		\label{eff_nonlin2 }
	    \end{subfigure}
	\caption{Computation of the coefficient $M(\xi)$ in Experiment set 3, with $\beta=0.1$, for $\tilde{\alpha}=0.05$ (a) and $\tilde{\alpha}=0.07$ (b), using a $L^1$-type regularizing term in the functional $J_{\alpha, 3}(M)$ with $\alpha = 0.001$. The absolute $L^2$-error results to be approximately 0.01 in both of the cases. }
	\label{effect_nonlinearity2}
\end{figure}

We further refine the approach by introducing the functional
\begin{align}\label{functional2c}
&J_{\alpha, 3}(M) := \frac12 \int_{L_0}^{L_1} | \eta(\xi, T, M) - \eta_T(\xi) |^2 d\xi +\frac12 \int_{L_0}^{L_1} | u(\xi, T, M) - u_T(\xi) |^2 d\xi + \\
&~~~~~~~~~~~~ \frac{\alpha}{2} \int_{L_0}^{L_1} | M(\xi) -1 | d\xi, \notag
\end{align}
where the regularizing term uses the $L^1$-norm and $\alpha = 0.001$. The numerical results, shown in Figure \ref{effect_nonlinearity2}, indicate even better accuracy. For $\tilde{\alpha}=0.05$ (315 iterations, $\Delta t = 20/1500$, $\Delta x = 60/500$) and for $\tilde{\alpha}=0.07$ (500 iterations, $\Delta t = 20/1700 \approx 0.0117$, $\Delta x = 60/700 \approx 0.0857$), the $L^2$-error is reduced to approximately 0.01 in both cases.

As previously discussed, the Tikhonov functionals considered here include a penalty on deviations from the constant component of the coefficient $M(\xi)$. This regularization serves to stabilize the reconstruction by reinforcing the coefficient around the reference value within the solution space. Without such regularization, the reconstructed coefficient may exhibit large deviations in its constant part, potentially increasing the computational effort and adversely affecting the convergence of the iterative scheme.

\begin{figure}[ht]
  \centering
		\includegraphics[width=10cm, height=8cm]{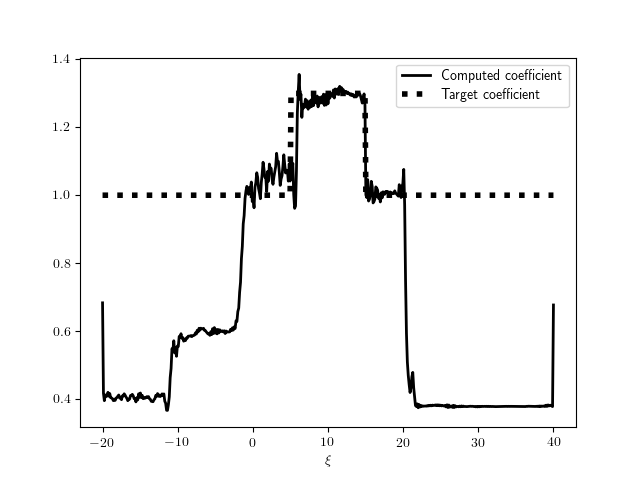}
	\caption{Computation of the coefficient $M(\xi)$ in Experiment set 3, obtained by minimizing the functional $J_{\alpha, 4}(M)$. }
	\label{Coeff_without_mean_regul}
\end{figure}

To illustrate this effect,  Figure \ref{Coeff_without_mean_regul} presents the results of a numerical experiment using the alternative functional
\begin{align}\label{functionalnomean}
&J_{\alpha, 4}(M) := \frac12 \int_{L_0}^{L_1} | \eta(\xi, T, M) - \eta_T(\xi) |^2 d\xi +\frac12 \int_{L_0}^{L_1} | u(\xi, T, M) - u_T(\xi) |^2 d\xi + \\
&~~~~~~~~~~~~ \frac{\alpha}{2} \int_{L_0}^{L_1} M(\xi)^2 d\xi, \notag
\end{align}
with a regularization parameter $\alpha = 0.001$. The experiment is carried out with $\tilde{\alpha} = 0.05$, $\beta = 0.1$, and the numerical parameters $\Delta t = 20/1500 \approx 0.013$ and $\Delta x = 60/500 = 0.12$. After 50 iterations, the recovered coefficient exhibits significant deviations, particularly in its constant part. The resulting absolute $L^2$-error is approximately 3.5, highlighting the importance of incorporating a regularization that reinforces the coefficient's alignment with the reference value.

It is important to note that the nonlinear regime extends beyond the theoretical results presented in this paper, which are derived for the linear case $\tilde{\alpha}= 0$. This suggests the potential extension of these analytical results to the full nonlinear, dispersive system \eqref{MA29}, which will be explored in a future research.


\vspace{0.5cm}

{\bf Experiment set 4:} Identification of other non-differentiable or discontinuous coefficients.

In this set of numerical experiments, we apply the coefficient reconstruction scheme to additional
 non-differentiable coefficient profiles. 
 
\begin{figure}
\centering
\includegraphics[width=10cm, height=8cm]{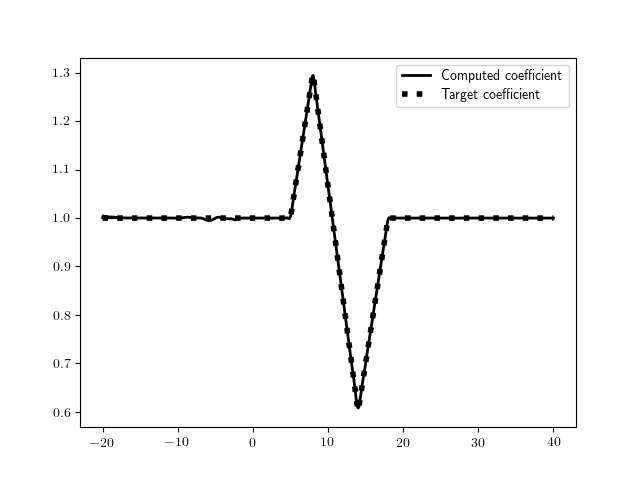}
\caption{Identification of a non-differentiable coefficient $M(\xi)$ in Experiment set 4, with parameters $\tilde{\alpha}=0, \beta=0.1$, and using functional $\mathcal{J}_{\alpha, 2}$ without regularization ($\alpha = 0$). }
\label{Coeff_it450_triangularcoeff_Error007}
\end{figure}

 We first consider the piecewise linear coefficient
\begin{equation}
M(\xi) = \begin{cases}
                 &1 + 0.1(\xi-5) , ~~~ 5\leq \xi \leq 8\\
                 &0.6 - (7/60)(\xi - 14), ~~~ 8 < \xi \leq 14\\
                 &1 + 0.1(\xi-18), ~~~ 14 < \xi \leq 18\\
                 &1, ~~~ \text{otherwise}.
             \end{cases}
\end{equation}
Figure \ref{Coeff_it450_triangularcoeff_Error007}, shows the result obtained after 450 iterations
using the parameters $\tilde{\alpha}=0$, $\beta = 0.1$, $\Delta t = 20/1700$, $\Delta x = 60/700$, and the functional $\mathcal{J}_{\alpha, 2}$ with $\alpha = 0$ (i.e. without regularization).
The reconstruction achieved an $L^2$-error of approximately 0.007.

\begin{figure}
\centering
\includegraphics[width=10cm, height=8cm]{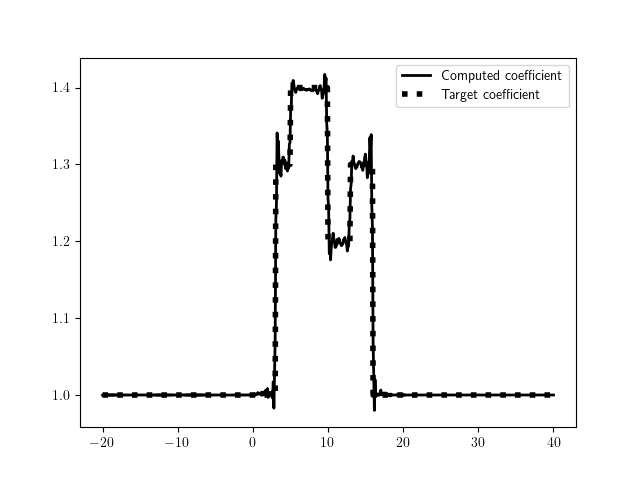}
\caption{Identification of a discontinuous coefficient $M(\xi)$ in Experiment set 4, with parameters $\tilde{\alpha}=0, \beta=0.1$, and using functional $\mathcal{J}_{\alpha,3}$ with $\alpha= 0.001$. }
\label{Coeff_it300_irregularcoeff_regul}
\end{figure}

Next, we consider the identification of a piecewise constant coefficient:
\begin{equation}
M(\xi) = \begin{cases}
                 &1.3 , ~~~ 3 \leq \xi \leq 5\\
                 &1.4, ~~~ 5 < \xi \leq 10\\
                 &1.2, ~~~ 10 < \xi \leq 13\\
                 &1.3, ~~~ 13 < \xi \leq 16\\
                 &1, ~~~ \text{otherwise}.
             \end{cases}
\end{equation}
The result of this reconstruction is displayed in Figure \ref{Coeff_it300_irregularcoeff_regul}, 
obtained after after 300 iterations using the same numerical parameters as in the previous experiment.
In this case, we used the Tikhonov functional $\mathcal{J}_{\alpha,3}$ with regularization parameter $\alpha=0.001$. The resulting $L^2$-error was approximately 0.06.

These results indicate that the proposed identification scheme performs well even for non-differentiable and discontinuous coefficients, highlighting the robustness of the numerical strategy developed in this work.

\vspace{0.5cm}

{\bf Experiment set 5:} Effect of noisy data on coefficient identification in the full nonlinear-dispersive system.

In this set of computer simulations, we investigate the impact of noisy measurements of wave amplitude and fluid velocity at the final time $t=T$, denoted as $\eta_T(\xi), u_T(\xi)$, respectively, in the context of the full nonlinear-dispersive Boussinesq system \eqref{MA29}. The aspect of noise in measurements is crucial in inverse problems, as errors invariably accompany measurements, stemming from unpredictable fluctuations in measurement devices or the experimenter's interpretation of instrumental readings. 

We set the modeling parameters $\tilde{\alpha}=\beta = 0.01$, and the numerical parameters $T=15$, $\Delta t = T/1500 = 0.01$, $\Delta x = 60/800 = 0.075$. The coefficient $M(\xi)$ corresponds to a chain of Gaussian profiles \eqref{gauss_coeff}. Additionally, Gaussian noise of magnitude $4 \times 10^{-2}$ is introduced along the interval $[-15, 30]$ and added to the observations $\eta_T(\xi)$ and $u_T(\xi)$.

\begin{figure}
\centering
\includegraphics[width=10cm, height=8cm]{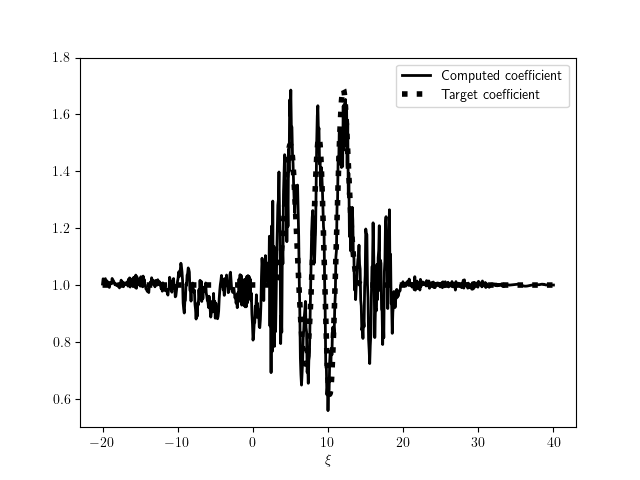}
\caption{Identification of coefficient $M(\xi)$ from noisy data in Experiment set 5, with parameters $\tilde{\alpha}=\beta=0.01$, and using functional \eqref{functional2a}, without regularization ($\alpha = 0$). }
\label{Exp4noiseNEW}
\end{figure}

Initially, coefficient recovery is performed by minimizing the regularized functional \eqref{functional2a} without a regularization term ($\alpha=0$), as in all previous experiments. The resulting coefficient after 20 iterations is depicted in Figure \ref{Exp4noiseNEW}, revealing a discrepancy of approximately 0.54 between the computed and exact coefficient $M(\xi)$.

\begin{figure}
\centering
\includegraphics[width=10cm, height=8cm]{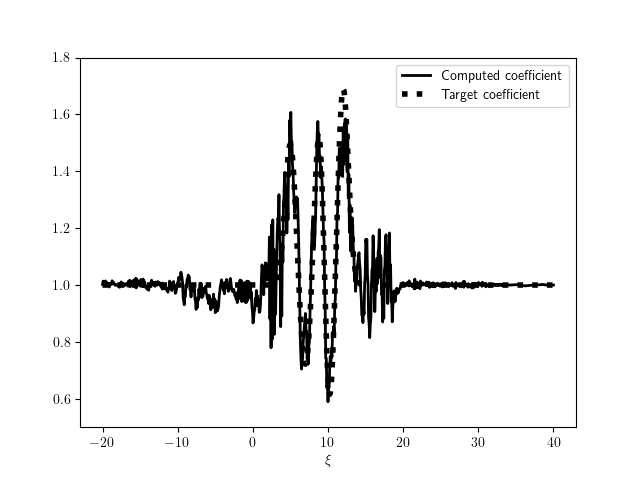}
\caption{Identification of coefficient $M(\xi)$ from noisy data in Experiment set 5, with parameters $\tilde{\alpha}=\beta=0.01$, and using functional \eqref{functional2a}, with regularization parameter $\alpha = 0.01$.}
\label{Exp6noiseNEW}
\end{figure}

Subsequently, we introduce a regularization term with $\alpha = 0.01$ in the functional $J_{\alpha,2}$, resulting in the coefficient shown in Figure \ref{Exp6noiseNEW} after 20 iterations of the L-BFGS-B scheme. The error between the computed and exact coefficients improves to approximately 0.43. This demonstrates that regularization in the Tikhonov functional \eqref{functional2a} indeed enhances the accuracy of the computed coefficient. 

\begin{figure}
\centering
\includegraphics[width=10cm, height=8cm]{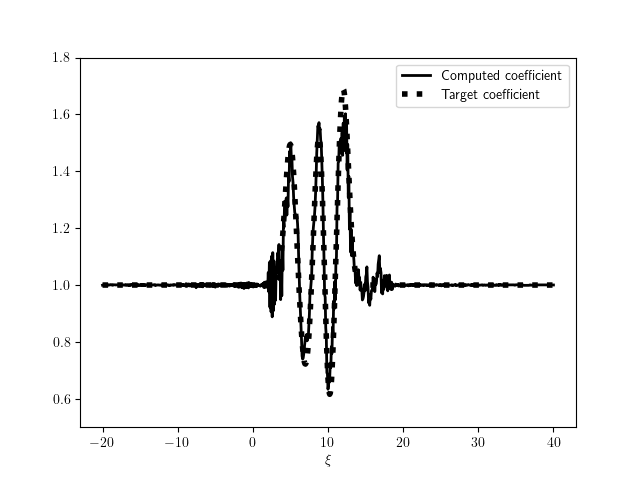}
\caption{Identification of coefficient $M(\xi)$ from noisy data in Experiment set 5, with parameters $\tilde{\alpha}=\beta=0.01$, and using functional \eqref{functional2c}, with regularization parameter $\alpha = 0.01$.}
\label{Exp5noiseNEW}
\end{figure}

Further enhancement is achieved by considering the functional \eqref{functional2c}
with $\alpha = 0.01$. Figure \ref{Exp5noiseNEW} shows the result of this simulation after 20 iterations, where
a significantly improved adjustment between the exact and computed coefficient is achieved 
with an absolute error of approximately 0.24, markedly lower than in previous numerical experiments.

From these experiments, it is evident that introducing a regularizing term is crucial, particularly in the presence
of noisy data where the ill-posed nature of the inverse problem is pronounced.

\section{Conclusions}

This paper investigated the problem of reconstructing the linear wave speed coefficient, $c(x)$, within a linear system Boussinesq-type model using measurements of water wave amplitude $N_T(\xi)$, and fluid velocity $V_T(\xi)$, taken at a final time $t=T$. We extended previous findings from \cite{JCPipicanoSosa} concerning the restricted case of a linear BBM equation. Our investigation transformed the inverse problem into a PDE-constrained optimization problem, focusing on an appropriate Tikhonov-type objective functional. This functional quantifies the disparity between measured data, $N_T(\xi), V_T(\xi)$, and actual data, $N(\xi, T; c), V(\xi, T; c)$.

Through analytical work, we established the existence of a minimum for the Tikhonov functional and demonstrated both local stability and uniqueness of solutions to the constrained optimization problem for fixed regularization parameter values. These results find further support through numerical illustrations, implemented using the FeniCS library to solve the boundary-initial value problem associated to the full nonlinear Boussinesq system \eqref{MA29}, and the L-BFGS-B routine for optimization, as facilitated by the Dolfin-Adjoint and Scipy libraries.

Additionally, we presented a series of numerical examples involving the full nonlinear Boussinesq system \eqref{MA29}, to illustrate the effectiveness, robustness, and accuracy of the proposed numerical strategy, even in scenarios where the unknown coefficient has discontinuities or when noise contaminates the final measurements. Notably, the techniques employed in this study in the linearized model \eqref{MA29} with $\tilde{\alpha}=0$, show promise for broader applications in coefficient recovery within other dispersive evolution equations. In particular, similar analytical approaches will be utilized to explore coefficient recovery for the fully nonlinear-dispersive Boussinesq system \eqref{MA29}, as planned in a forthcoming work.


\section*{Acknowledgments}
This work was partially supported by Universidad del Valle, Cali, Colombia, under research projects C.I. 71235, 71288, Cali, Colombia, and MinCiencias under project FP44842-266-2017.

\end{document}